\pdfoutput=1
\RequirePackage{ifpdf}
\ifpdf % We~are running pdfTeX in pdf mode
\documentclass[pdftex]{sigma}
\else
\documentclass{sigma}
\fi

\numberwithin{equation}{section}

\newtheorem{Theorem}{Theorem}[section]
\newtheorem*{Theorem*}{Theorem}

\newtheorem{Lemma}[Theorem]{Lemma}
\newtheorem{Proposition}[Theorem]{Proposition}

\theoremstyle{definition}

\newtheorem{Remark}[Theorem]{Remark}

\usepackage{mathrsfs}

\begin{document}
%\allowdisplaybreaks

\newcommand{\arXivNumber}{2502.03351}

\renewcommand{\PaperNumber}{018}

\FirstPageHeading

\ShortArticleName{Regularized $\zeta_{\Delta}(1)$ for Polyhedra}

\ArticleName{Regularized $\boldsymbol{\zeta_{\Delta}(1)}$ for Polyhedra}

\Author{Alexey Yu. KOKOTOV and Dmitrii V. KORIKOV}
\AuthorNameForHeading{A.Yu.~Kokotov and D.V.~Korikov}
\Address{Department of Mathematics and Statistics, Concordia University,\\
1400 De Maisonneuve Blvd.~W., Montreal, QC H3G 1M8, Canada}
\Email{\mail{alexey.kokotov@concordia.ca}, \mail{dmitrii.v.korikov@gmail.com}}

\ArticleDates{Received May 13, 2025, in final form February 03, 2026; Published online February 25, 2026}

\Abstract{Let $X$ be a compact polyhedral surface (a compact Riemann surface with flat conformal metric $\mathfrak{T}$ having conical singularities). The $\zeta$-function $\zeta_\Delta(s)$ of the Friedrichs Laplacian on $X$ is meromorphic in ${\mathbb C}$ with a single simple pole at $s=1$. We define~$\operatorname{reg}\zeta_\Delta(1)$ as \smash{$\lim\limits_{s\to 1} \bigl( \zeta_\Delta(s)-\frac{ {\rm Area}(X,\mathfrak{T}) }{4\pi(s-1)}\bigr)$}. We derive an explicit expression for this spectral invariant through the holomorphic invariants of the Riemann surface $X$ and the (generalized) divisor of the conical points of the metric $\mathfrak{T}$. We study the asymptotics of $\operatorname{reg}\zeta_\Delta(1)$ for the polyhedron obtained by sewing two other polyhedra along segments of small length. In addition, we calculate $\operatorname{reg}\zeta(1)$ for a family of (non-Friedrichs) self-adjoint extensions of the Laplacian on the tetrahedron with all the conical angles equal to $\pi$.}


\Classification{58J52; 35P99; 30F10; 30F45; 32G15; 32G08}

\section{Introduction}

Let $X$ be a compact Riemann surface of genus $g$ and let $\sum_{k=1}^Mb_k P_k$ ($P_k\in X$, $b_k\in {\mathbb R}$, ${b_k>-1}$, $M$ is a positive integer) be a (generalized) divisor of degree $\sum_{k=1}^M b_k=2g-2$. Then (see~\cite{Troyanov})\looseness=1
\begin{itemize}\itemsep=0pt
\item there exists a unique (up to rescaling) conformal metric $\mathfrak{T}$ (the Troyanov metric in the sequel) on $X$ such that $\mathfrak{T}$ is flat
in $X\setminus\{P_1, \dots, P_M\}$ and $\mathfrak{T}$ has conical points at $P_k$ of conical angles $\beta_k=2\pi(b_k+1) $. There exist holomorphic coordinates $x_k$ near $P_k$ (so called distinguished local parameters on $X$) such that \smash{$\mathfrak{T}=\big|{\rm d}\bigl(x_k^{\beta_k/2\pi}\bigr)\big|^2=(b_k+1)^2|x_k|^{2b_k}|{\rm d}x_k|^2$}.

\item there exists a triangulation of $X$ consisting of Euclidean triangles such that all the $P_k$ are among the vertices of the triangles from this triangulation.
\end{itemize}

Clearly, any polyhedral surface (i.e., a closed orientable compact surface glued from Euclidean triangles) can be provided with natural complex structure and flat conformal conical metric. This results in a pair $(X, \mathfrak{T})$ described above.

Let $\dot{X}:=X\setminus \{P_1, \dots, P_M\}$ and let $\Delta$ be the Friedrichs self-adjoint extension of the symmetric Laplacian $-4\mathfrak{T}^{-1}\partial\overline{\partial}=:\Delta_0$ on $X$ corresponding to a Troyanov metric $\mathfrak{T}$ with domain $C^\infty_0(\dot{X})$. Then (see, e.g., \cite{KokLN})
\begin{itemize}\itemsep=0pt
\item $\Delta$ has a discrete spectrum consisting of eigenvalues, $\lambda_k$, of finite multiplicity,
\item the operator $\zeta$-function
\begin{equation}
\label{zeta}
\zeta_\Delta(s)=\sum_{\lambda_k>0}\frac{1}{\lambda_k^s}
\end{equation}
defined initially for $\operatorname{Re} s>1$ can be analytically continued to a meromorphic function in~${\mathbb C}$ with unique (simple) pole at $s=1$. One has the representation
\begin{equation*}
\zeta_\Delta(s)=\frac{1}{\Gamma(s)}\left\{\frac{A}{4\pi(s-1)}+
\left[\frac{1}{12}\sum_{k=1}^M\left\{\frac{2\pi}{\beta_k}-\frac{\beta_k}{2\pi}\right\}-1 \right]\frac{1}{s}+e(s) \right\},
\end{equation*}
where $A:={\rm Area }(X,\mathfrak{T})$ is the area of the surface $X$ in the metric $\mathfrak{T}$ and $e(s)$ is an entire function.
\end{itemize}
In particular, $\zeta_\Delta(s)$ is regular at $s=0$, this gives an opportunity to introduce the spectral determinant of the Laplacian on a polyhedral surface, $\det \Delta$, as $\exp{(-\zeta_\Delta'(0))}$. This quantity was explicitly computed and extensively studied in \cite{KokTAMS,KokPAMS,KokKorJDG}. The main goal of the present paper is to study another interesting spectral invariant of the polyhedral surface $X$, namely the regularized value of the operator $\zeta$-function
at $s=1$
\begin{equation}
\label{reg zeta of 1}
\operatorname{reg}\zeta_\Delta(1)= \lim_{s\to 1} \left( \zeta_\Delta(s)-\frac{A}{4\pi(s-1)}\right) .
\end{equation}

\begin{Remark}
We choose regularization \eqref{reg zeta of 1}, following the tradition (in the case of smooth metrics this regularization was previously used in \cite{Okikitorus,Okikiolu,Steiner}); however, it is worth noting that a slightly different regularization $\widetilde{\operatorname{reg}}\,\zeta_\Delta(1)=\operatorname{reg}\zeta_\Delta(1)-A\log A/4\pi$ (obtained by replacement~$A$ with $A^s$ in \eqref{reg zeta of 1}) satisfies the same scaling law $\widetilde{\operatorname{reg}}\,\zeta_{\Delta'}(s)=t\, \widetilde{\operatorname{reg}}\,\zeta_{\Delta}(s)$ of $\widetilde{\operatorname{reg}}\,\zeta_\Delta(1)$ under the homothety $\mathfrak{T}':=t\mathfrak{T}$ as the $\zeta$-function itself $\zeta_{\Delta'}(s)=t^s\zeta_{\Delta}(s)$.
\end{Remark}

It should be noted that for the classical Riemann $\zeta$-function, $\zeta_R(s)$, one has
\[\operatorname{reg}\zeta_R(1)=\gamma=0.577 \dots . \]
So in a sense we are dealing with an analogue of the Euler--Mascheroni constant for a polyhedral surface.

We derive an explicit expression for $\operatorname{reg}\zeta_\Delta(1)$ through the holomorphic invariants of $X$ and the (generalized) divisor of the conical points of the Troyanov metric $\mathfrak{T}$. The derivation consists of the following three ingredients.

{\bf 1.} Let $G$ denote the Green function for the Laplacian $\Delta$ (see, e.g., \cite[p.~25]{Fay}) and let $G_s$ ($s=2,3,\dots$) denote the Green functions for $\Delta^s$. Recall that the values of $\zeta_\Delta(s)$ at $s=2,3,\dots$ are connected with diagonal values $G_s(x,x)$ via
\[\zeta_\Delta(s)={\rm Tr}\Delta^{-s}=\int_{X}G_s(x,x)\,{\rm d }S(x), \qquad s=2,3,\dots.\]
(Here and in what follows ${\rm d }S(x)$ (sometimes ${\rm d }S$) denotes the area element.)

A version of the above formula for $s=1$ is derived by Steiner \cite{Steiner} for the case of {\it smooth} metrics on $X$. In this case, both $\zeta_\Delta(1)$ and $G(x,x)$ are replaced by their regularized values. Namely, introduce the {\it Robin mass}
\begin{equation}
\label{Robin mass}
m(P)=\lim_{P'\to P}\left[G(P,P')+\frac{\log d(P,P')}{2\pi}\right], \qquad P\in X,
\end{equation}
where $d(P,P')$ is the geodesic distance between $P$ and $P'$. Then (see \cite[Proposition 2]{Steiner})
\begin{equation}
\label{Zeta Robin}
\operatorname{reg}\zeta_\Delta(1)=\int_{X}\left(m-\frac{\log 2-\gamma}{2\pi}\right){\rm d }S .
\end{equation}
 Steiner's proof of \eqref{Zeta Robin} is based on the classical short-time asymptotics of the heat kernel $H_t$ corresponding to a smooth metric.

In Section \ref{Steiner sec}, we prove the first main result of the present paper.
\begin{Theorem}
Steiner's relation \eqref{Zeta Robin} holds for the case of the $($singular$)$ Troyanov metric $\mathfrak{T}$.
\end{Theorem}

 To implement (in a more sophisticated way) Steiner's argument, we derive the corrections to the parametrix of $H_t$ arising due to the presence of conical singularities of the metric. To this end, we make use of the explicit solutions to the heat equation in cones \cite{Carslaw,Dowker}.

Due to \eqref{Zeta Robin}, the calculation of $\operatorname{reg}\zeta_\Delta(1)$ is reduced to the calculation of the Robin mass $m$.

{\bf 2.} The Robin mass $m$ is calculated by the use of the following Verlinde--Verlinde (simply Verlinde in the sequel) formula \eqref{Green Verlinde} (see \cite[equations~(5.7)--(5.10)]{Ver}; it should be noted that usually this name is associated to the famous formula of the same authors from conformal field theory that has nothing to do with \eqref{Green Verlinde}).

Let $\rho^{-2}(z)|{\rm d}z|^2$ be a conformal metric on $X$ (the case of the Troyanov metric $\mathfrak{T}$ is also allowed, see Lemma \ref{Ver Vert Lem} below); introduce the prime form $E(x,y)$, the Abel transform $\mathscr{A}$ and the matrix of $b$-periods $\mathbb{B}$ for the (marked) Riemann surface $X$. Introduce a symmetric section
\begin{equation}
\label{Verlinde}
F(x,y)=\exp\bigl[-2\pi\operatorname{Im}\mathscr{A}(x-y)^{\mathsf T} (\operatorname{Im}\mathbb{B})^{-1}\operatorname{Im}\mathscr{A}(x-y)\bigr]|E(x,y)|^2
\end{equation}
of $|K_x|^{-1}\hat{\otimes}|K_y|^{-1}$ (here $K$ denotes the canonical bundle; the notation means that $F$ is a section of the bundle $|K|^{-1}$ with respect to both arguments $x$ and $y$, cf.\ \cite[p.\ 23]{Fay}). Then (see the last paragraph of \cite[Section~5]{Ver})
\begin{equation}\label{Green Verlinde}
G(x,y)=\frac{1}{2}\bigl(m(x)+m(y)\bigr)+\Phi(x,y),
\end{equation}
where
\begin{equation}\label{Phi func}
\Phi(x,y)=-\frac{1}{4\pi}\log \frac{F(x,y)}{\rho(x)\rho(y)}.
\end{equation}
Integrating \eqref{Green Verlinde} over $x$ and $y$ and taking into account that $G(\cdot,y)$ is $L_2$-orthogonal to constants, we obtain
\begin{align}\label{scalar Robin 2}
\int_{X}m\,{\rm d }S=&-\frac{1}{A}\int_{X}\int_{X}\Phi(x,y)\, {\rm d }S(x){\rm d }S(y),
\end{align}
where $A$ is the area of $X$ in the metric $\rho^{-2}(z)|{\rm d}z|^2$. Using \eqref{scalar Robin 2}, one derives the expressions for the Robin mass $m$ and $\operatorname{reg}\zeta_\Delta(1)$ via the metric on $X$. Namely, integrating \eqref{Green Verlinde} with respect to~$y$, one gets
\[%\label{samaMassa}
m(x)=\frac{1}{A^2}\int_X \int_X\Phi(x, y) \,{\rm d }S(x) {\rm d }S(y)-\frac{2}{A}\int_X\Phi(x, y)\, {\rm d }S(y),
\]
whereas
combining \eqref{Zeta Robin} with \eqref{scalar Robin 2} and \eqref{Phi func}, \eqref{Verlinde} yields
the following result.
\begin{Theorem}
Let $\rho^{-2}(z)|{\rm d}z|^2$ be either a smooth conformal metric or a Troyanov metric on~$X$.
One has the following explicit expression for the regularized value of the $\zeta_\Delta(s)$ at $s=1$:
\begin{align}
\operatorname{reg}\zeta_\Delta(1)={}&\frac{A}{2\pi}\left[\gamma-\log 2+\int_{X}\int_{X}\left(\log \frac{|E(x,y)|}{\rho^{1/2}(x)\rho^{1/2}(y)}\right.\right.\nonumber\\
&\left.\left.-2\pi\operatorname{Im}\mathscr{A}(x-y)^{\mathsf T}(\operatorname{Im}\mathbb{B})^{-1}\operatorname{Im}\mathscr{A}(x-y)\right)\frac{{\rm d }S(x) {\rm d }S(y)}{A^2}\right] .
\label{zeta of 1 via metric}
\end{align}
\end{Theorem}

{\bf 3.} An explicit formula for the Troyanov metric (in higher genus) was first proposed in \cite{KokLN}; in~\cite{KokPAMS}, it was given a more effective one but only for a metric with no Weierstrass points among the vertices $P_1,\dots,P_M$. In Section \ref{Troyanov sec}, we improve the construction from \cite{KokPAMS}: in contrast to \cite{KokPAMS} the new expression for $\mathfrak{T}$ given by formulas \eqref{Troy 1}, \eqref{Troy 2}, \eqref{Troy 3} below is valid in general. We also recall formulas \eqref{Troy sphere} and \eqref{Troy tori} for the case $g=0$ and $g=1$, respectively. The explicit expression for $\operatorname{reg}\zeta_\Delta(1)$ for polyhedra is obtained by substituting these formulas into \eqref{zeta of 1 via metric}.

Next, in Section \ref{asymp sec}, we consider a family $\varepsilon\mapsto X(\varepsilon)$ of polyhedral surfaces $X(\varepsilon)$ with two conical points $P(\varepsilon)$ and $P'(\varepsilon)$ (with conical angles $4\pi$) approaching to each other in such a way that $\varepsilon={\rm dist}(P(\varepsilon),P'(\varepsilon))$ and the ``limit surface'' $X(0)$ splits into two disjoint surfaces $X_\pm$ with punctures at $P_\pm$, respectively. We prove the following result.
\begin{Theorem} One has the asymptotic formula
\begin{equation}
\label{zeta 1 asymp}
\lim_{\varepsilon\to 0}\left[\operatorname{reg} \zeta_{\Delta}(1)+\frac{A_+ A_-}{A}\left(\frac{\log (\varepsilon/4)}{\pi}-\sum_\pm m_\pm(P_\pm)\right)\right]=\sum_{\pm}\operatorname{reg} \zeta_{\Delta_{\pm}}(1)
\end{equation}
for the regularized value at $s=1$ of the $\zeta$-function of the Laplacian $\Delta:=\Delta(\varepsilon)$ on $X:=X(\varepsilon)$. Here $A$ and $A_\pm$ are areas of $X$ and $X_\pm$, respectively, $m_\pm(P_\pm)$ are the Robin masses at $P_\pm\in X_\pm$.
\end{Theorem}

We also derive (via the matched asymptotics method, see \cite[Chapters 2 and 6]{MNP}) the asymptotics of the first non-zero eigenvalue of $\Delta(\varepsilon)$,
\begin{equation}
\label{first eigenvalue asym}
\frac{1}{\lambda_1(\varepsilon)}=-\frac{A_+A_-}{A}\left[\frac{1}{\pi}\log (\varepsilon/4)-\sum_\pm m_\pm(P_\pm)\right]+o(1).
\end{equation}
Using this asymptotics, one observes that \eqref{zeta 1 asymp} takes the form
\begin{equation}
\label{zeta_1 asymp_1}
\lim_{\varepsilon\to 0}\bigl[\operatorname{reg} \zeta_{\Delta}(1)-\lambda_1(\varepsilon)^{-1}\bigr]=\sum_{\pm}\operatorname{reg} \zeta_{\Delta_{\pm}}(1),
\end{equation}
where $\Delta_\pm$ are the Friedrichs Laplacians on $X_\pm$, respectively. Formula \eqref{zeta_1 asymp_1} exhibits additive law for $\operatorname{reg} \zeta_{\Delta}(1)$ after the exclusion of not only the zero mode but also the mode with eigenvalue~$\lambda_1(\varepsilon)$ tending to zero as $\varepsilon\to 0$. The reason for this additional exclusion is that in the right-hand side of \eqref{zeta_1 asymp_1} the zero modes are excluded from both summands $\operatorname{reg} \zeta_{\Delta_{\pm}}(1)$. Although \eqref{zeta_1 asymp_1} follows from the formal substitution of the asymptotics
\begin{equation}
\label{separated eigenvalues asym}
\lambda_{k+1}(\varepsilon)=\lambda_{k}+o(1)
\end{equation}
(obtained again by the matched asymptotics method; here $k>1$ and $\lambda_1\le\lambda_2\le\cdots$ is the non-decreasing sequence of non-zero eigenvalues of the operator ${\rm diag}(\Delta_+,\Delta_-)$ in $L_2(X_+)\times L_2(X_-)$) into \eqref{zeta}, it cannot be justified in such a way since asymptotics \eqref{separated eigenvalues asym} are not uniform in $k$. Moreover, such a substitution is incorrect in general: indeed, it would lead to the incorrect (inverse logarithmic due to \eqref{first eigenvalue asym} instead of power-type) dependence on the parameter $\varepsilon$ in \cite[formula (3.15)]{KokTAMS} for the determinant $\det \Delta(\varepsilon)$.

Finally, in Section \ref{Examples sec}, we consider the family $\alpha\mapsto\Delta_{4,\alpha}$ ($\alpha\in[0,\pi)$) of self-adjoint extensions of the Laplacian on the tetrahedron $\mathscr{T}=ABCD$ with domain consisting of functions on $\mathscr{T}\simeq \mathbb{C}P^1$ vanishing near one of the vertices (say, $A$). We restrict ourselves to the case in which all the conical angles of $\mathscr{T}$ are equal to $\pi$; in this case, $\mathscr{T}$ is isometric to the Riemann sphere equipped with metric $C(\mathscr{T})|\mathfrak{Q}|$, where $\mathfrak{Q}$ is the quadratic differential
\[\mathfrak{Q}=\frac{({\rm d}z)^2}{z(z-1)(z-\xi(\mathscr{T}))}\]
with $\xi(\mathscr{T})\in\mathbb{C}\backslash\{0,1\}$; $C(\mathscr{T})>0$.
 The domain of $\Delta_{4,\alpha}$ contains all the smooth (in $\mathscr{T}\setminus \{A\}$; the
 $C^\infty$-structure on $\mathscr{T}$ comes from the Riemann sphere) functions with asymptotics{\samepage
\[u(x)={\rm const}\left(\frac{\log d(x,A)}{2\pi}\sin(\alpha)+\cos (\alpha)\right)+o(1)\]
as $x\to A$. The Friedrichs Laplacian, $\Delta_4$, coincides with $\Delta_{4,0}$.}

For the above tetrahedra $ABCD$ (with all four conical angles equal to $\pi$) having the unit area and the aspect ratio
\begin{equation}
\label{aspect ratio}
AB:BC:CA=|\tau|:|\tau-1|:1
\end{equation}
(where $\operatorname{Im}\tau>0$), we prove the following statement.
\begin{Theorem}
The explicit formulas for the regularized $\zeta$-functions of $\Delta_{4,\alpha}$ hold
\begin{gather}
\label{zeta tetra}
\operatorname{reg}\zeta_{\Delta_4}(1)=\frac{1}{\pi}\left[\frac{\gamma}{2}-\frac{\log (2\operatorname{Im}\tau)}{4}-\log |2\pi\eta(\tau)|\right],\\
\label{reg zeta for pseudoL on tetra}
\operatorname{reg}\zeta_{\Delta_{4,\alpha}}(1)=\frac{1}{\pi}\left[\frac{\gamma}{2}-\frac{\log (2\operatorname{Im}\tau)}{4}+\log |2\pi\eta(\tau)|\right]-2\operatorname{cot}\alpha, \qquad \alpha\ne 0 ,
\end{gather}
where $\eta$ is the Dedekind eta function.
\end{Theorem}

Note that $\operatorname{reg}\zeta_{\Delta_{4,\alpha}}(1)$ of $\Delta_{4,\alpha}$ tends to infinity as $\alpha\to +0$ since the first two eigenvalues of $\Delta_{4,\alpha}$ obey $\lambda_{1}(\Delta_{4,\alpha})\to -0$ (mind the minus sign: the first eigenvalue of $\Delta_{4, \alpha}$ is negative for~${\alpha>0}$), $\lambda_{2}(\Delta_{4,\alpha})\to +0$ as $\alpha\to+0$ (see \cite[Section 2]{CdV}) while the zero mode is excluded from~$\operatorname{reg}\zeta_{\Delta_4}(1)$.

We also mention in passing that there are conclusive numerical arguments in favour of the statement that among all the above tetrahedra, the (global) minimum of $\operatorname{reg}\zeta_{\Delta_4}(1)$ is attained at the regular one.

We hope that the results of the present paper can find their applications to extremal geometries of polyhedral surfaces (see \cite{Sarnak}). Special values of the $\zeta$-function are among the most interesting functionals on various moduli spaces and their extremal properties are intimately related to geometric properties of these spaces. Special values of the $\zeta$-function may be applied in the studies of inverse spectral problems in the spirit of Osgood--Phillips--Sarnak \cite{Zeld,OPS}. As~an additional motivation, we notice that in \cite{KokKorik}, using the Ricci flow, it was shown that a genus zero version of Theorem 1.3 implies the old result of Morpurgo~\cite{Mor} concerning the extremality of $\operatorname{reg} \zeta_\Delta(1)$ at the round metric on the Riemann sphere.

\section{Explicit expressions for Troyanov metric}
\label{Troyanov sec}
In this section, we derive the explicit expressions for the Troyanov metric $\mathfrak{T}$ via its divisor.

\subsection*{The case $\boldsymbol{g>1}$} Choose a marking on $X$ (i.e., a base point $p_0$ and $2g$ loops from
$\pi_1(X, p_0)$ forming a canonical homotopy basis; these loops can be chosen in such a way that dissecting $X$ along them yields a fundamental polygon $\mathcal{D}$; see \cite{Fay}). In order to construct $\mathfrak{T}$, it suffices to construct, for each~${k=1,\dots,M}$, a~multiplicative holomorphic differential of weight $1$, $\mathcal{P}(\cdot, P_k)$, on $X$ with a~single zero of multiplicity $2g-2$ at $P_k$ and {\it unitary} multipliers along the basic cycles.
Then the Troyanov metric is given by
\begin{equation}
\label{Troy 1}
\mathfrak{T}=C \prod_{k=1}^M|\mathcal{P}(\cdot, P_k)|^{\frac{2b_k}{2g-2}}\end{equation}
with arbitrary positive constant $C$ responsible for the rescaling of the metric.

To construct the $\mathcal{P}(z, P)$ ($P\in X$) introduce (following \cite[formula (1.13)]{Fay}, see also \cite{Ver})
\begin{equation}
\label{Troy 2}
\sigma(z,z')=\frac{\theta\bigl(\mathscr{A}\bigl(\sum_{j=1}^g z_j-gz\bigr)+K^z\bigr)}{\theta\bigl(\mathscr{A}\bigl(\sum_{j=1}^g z_j-z'-(g-1)z\bigr)+K^z\bigr)}\prod_{j=1}^g\frac{E(z_j,z')}{E(z_j,z)},
\end{equation}
where $z_1,\dots,z_g,z,z'\in\mathcal{D}$, $\mathscr{A}$ and $K^{z}$ are the Abel transform (with integration going along the paths lying inside $\mathcal{D}$) and the vector of Riemann constants corresponding to a base point $z\in \mathcal{D}$; $E(\cdot, \cdot)$ is the prime form; $\theta$ denotes the theta function of $X$ (with respect to the above marking).

It can be checked that $\sigma(z, z')$ is independent of $z_1,\dots,z_g$ and presents a nowhere vanishing multiplicative holomorphic $\frac{g}{2}$-differential with respect to the first argument, $z$, with multipliers~$1$, $ \exp((g-1)\pi {\rm i} \mathbb{B}_{jj}+2\pi {\rm i} K^z_j)$
along the cycles $a_j$, $b_j$, respectively (see \cite[Proposition 1.2]{Fay}). Here~$\mathbb{B}$ is the matrix of $b$-periods of $X$.

Let
\begin{gather}
\label{VerFay}
\mathcal{P}(z,z'(P))=\sigma^2(z,z_0)E^{2g-2}(z,z'(P))\exp\left(-4\pi {\rm i} \sum_{s=1}^{g}\rho_s\int_{z_0}^z v_s\right), \end{gather}
where $z_0$ is an arbitrary point of $X$ (the value of $\sigma$ is computed with respect to an arbitrary holomorphic local parameter at $z_0$), and $\vec{\rho}=(\rho_1,\dots,\rho_g)^{\mathsf T}\in\mathbb{R}^g$; $\{v_s\}$ is the canonical basis of Abelian differentials (corresponding to the chosen marking of $X$).

(Notice that the expression $\sigma(z,z_0)E^{g-1}(z,z'(P))$ appeared previously in \cite{Ver} and, later, in~\cite{Fay}. We introduce the exponential factor in the right hand side of \eqref{VerFay} to make the monodromies unitary.)

Recall that the prime form $E(z,z'(P))$ is a multiplicative holomorphic $(-1/2)$-differential on~$X$ with the unique zero $P$ of order $1$ and with multipliers
\[1, \qquad \exp\left(-\pi {\rm i} \mathbb{B}_{jj}-2\pi {\rm i} \int_{z'}^z v_j\right)\]
along the cycles $a_j$, $b_j$, respectively (see \cite[equation~(1.4)]{Fay}). Then $\mathcal{P}(\cdot,z'(P))$ is a multiplicative holomorphic $1$-differential on $X$ with the unique zero $P$ of order $2g-2$ and with multipliers
\begin{equation}
\label{multipliers}
\exp(-4\pi {\rm i} \rho_j), \qquad \exp\left(4\pi {\rm i} \left(K_j^{P}-\sum_{s=1}^{g}\rho_s\mathbb{B}_{sj}\right)\right)
\end{equation}
along the cycles $a_j$, $b_j$, respectively. To check this one uses the formula
\[K^z-(g-1)\int_{z'}^z \vec{v}=K^{z'}\]
(see \cite[Proposition 1.1]{Fay}). The multipliers \eqref{multipliers} are unitary if and only if
$\operatorname{Im}\mathbb{B}\vec{\rho}=\operatorname{Im} K^{P}$.
Therefore, the expression
\begin{equation}
\label{Troy 3}
\mathcal{P}(z , P)=\sigma^2(z ,z_0)E^{2g-2}(z ,z'(P))\exp\bigl(-4\pi {\rm i} \mathscr{A}(z -z_0)^{\mathsf T}(\operatorname{Im}\mathbb{B})^{-1}\operatorname{Im} K^{P}\bigr)
\end{equation}
(where the holomorphic local coordinate $z'$ at $P$ is arbitrary) gives a multivalued holomorphic (Prym) differential on $X$ with the unique zero $P$ of order $2g-2$ and with unitary multiplicative twists along the basic cycles.

\subsection*{The case $\boldsymbol{g=1}$} In this case, the construction is similar (and was done in \cite[Section 4.2.3]{PirioFlatTori}). The idea is essentially the same: in genus $1$ one makes use of the theta function \smash{$\theta\bigl[^{1/2}_{1/2}\bigr](\cdot -P)$} with odd characteristics (i.e., the main ingredient of the prime form on the elliptic curve $X$ that vanishes at $P$) and introduces some exponential factor to make the twists unitary.

Namely, the product
\[\prod_{k=1}^M \left( \theta\left[^{1/2}_{1/2}\right](z-P_k)\right)^{b_k}\]
has a unitary multiplicative twist along the $a$-cycle and
the multiplicative twist
\begin{align*}
\exp{\left(-\pi {\rm i} \sum_{k=1}^Mb_k\bigl({\mathbb B}+2(z-P_k)+1\bigr)\right)}&=\exp{\left(2\pi {\rm i} \sum_{k=1}^Mb_k P_k\right)}\\
&=\exp{({\rm i}\psi)}\exp{\left(-2\pi \sum_{k=1}^Mb_k \operatorname{Im} P_k\right)}
\end{align*}
along the $b$-cycle (here ${\mathbb B}$ is the $b$-period of the elliptic curve $X$, we used the equality $\sum_{k=1}^M b_k=0$) with some real $\psi$.

For $\alpha\in {\mathbb R}$, the function $\exp{({\rm i}\alpha z)}$ has a unitary multiplicative twist along the $a$-cycle ($z\mapsto z+1$) and the multiplicative twist $\exp{({\rm i}\alpha {\mathbb B})}$ along the $b$-cycle
($z\mapsto z+{\mathbb B}$). Choose
\[\alpha=-\frac{2\pi}{\operatorname{Im}{\mathbb B}}\sum_{k=1}^Mb_k\operatorname{Im} P_k.\]
Then the Troyanov metric is given by
\begin{equation}
\label{Troy tori}
\mathfrak{T}=C\left|{\rm e}^{{\rm i}\alpha z} \prod_{k=1}^M \left\{ \theta\left[^{1/2}_{1/2}\right](z-P_k)\right\}^{b_k} {\rm d}z \right|^2 .
\end{equation}

\subsection*{The case $\boldsymbol{g=0}$} In genus zero such a formula is obvious (and was given in \cite{Troyanov})
\begin{equation}
\label{Troy sphere}
\mathfrak{T}=C\prod_{k=1}^M \left|z-z(P_k)\right|^{2b_k}|{\rm d}z|^2 ,
\end{equation}
where $z$ is the global coordinate on ${\mathbb C}\subset {\mathbb C}P^1=X$.

\subsection*{Verlinde formula (\ref{Green Verlinde}) for polyhedra} In this subsection, we prove the following statement.
\begin{Lemma}\label{Ver Vert Lem}
Formula {\rm \eqref{Green Verlinde}} is valid if the metric on $X$ is Troyanov.
\end{Lemma}
\begin{proof}
Let $y\in\dot{X}$. In view of \eqref{Phi func} and \eqref{Verlinde}, we have
\begin{equation}
\label{delta of phi}
\Delta_{x}\bigl(G(x,y)-\Phi(x,y)\bigr)=\mathfrak{T}^{-1}(x)B(x,x)-A^{-1}, \qquad x\in\dot{X},
\end{equation}
where $\vec{v}=(v_1,\dots,v_g)^{\mathsf T}$ and
\smash{$
%\label{Bergman}
B(x,y)=\overline{\vec{v}(y)}^{\mathsf T}(\operatorname{Im}\mathbb{B})^{-1}\vec{v}(x)
$}
is the Bergman reproducing kernel for holomorphic differentials.

Since the Green function $G$ corresponds to the Friedrichs Laplacian, $G(\cdot,y)$ is bounded near vertices. Therefore, \eqref{Phi func} implies the asymptotics
\begin{align}
G(x_k,P_k)-\Phi(x_k,P_k)&=\frac{1}{8\pi}\log \mathfrak{T}(x_k)+O(1)\nonumber\\
&=\frac{b_k}{4\pi}\log |x_k|+O(1)=\frac{b_k\log d(x_k,P_k)}{4\pi(b_k+1)}+O(1)\label{Phi near vertices}
\end{align}
near $P_k$ ($k=1,\dots,M$). Denote
\[\mathfrak{D}(x)=\sum_k \frac{b_k}{4\pi(b_k+1)}\chi_k(x)\log d(x,P_k),\]
where $\chi_k$ is a smooth cut-off function equal to 1 near $P_k$ and such that the diameter $\operatorname{supp}\chi_k$ (in the metric $\mathfrak{T}$) is sufficiently small. Then $\mathfrak{F}:=\Delta \mathfrak{D}\in C_c^{\infty}(\dot{X})$. It is easy to check that
\begin{equation}\label{for ortho}
\int_{X}\mathfrak{F}(x)\,{\rm d }S(x)=\sum_k \frac{b_k}{4\pi(b_k+1)}\beta_k=\frac{1}{2}\sum_{k}b_k=g-1.
\end{equation}

Now, introduce a solution $\tilde{m}$ to
\begin{equation}\label{correction robin eq}
\Delta\tilde{m}(z)=\mathfrak{T}^{-1}B(z,z)-A^{-1}-\mathfrak{F}(z),
\end{equation}
where $\Delta$ is the Friedrichs Laplacian. From \eqref{for ortho} and the well-known equality (see, e.g., \cite[p.15]{Fay})
$\int_X B(z,z)=g$,
it follows that the right-hand side of \eqref{correction robin eq} is $L_2(X,\mathfrak{T})$-orthogonal to constants. Thus, a solution to \eqref{correction robin eq} exists and is unique up to an additive constant. Denote $m=2(\mathfrak{D}+\tilde{m})$ and put $\mathfrak{m}(x,y):=G(x,y)-\Phi(x,y)-\frac{1}{2}(m(x)+m(y))$. Then
$
\Delta_x\mathfrak{m}(x,y)=\Delta_y\mathfrak{m}(x,y)=0$, $x,y\in\dot{X}$
and $\mathfrak{m}(\cdot,y)$, $\mathfrak{m}(x,\cdot)$ are bounded near vertices. Therefore, $\mathfrak{m}(\cdot,y)$, $\mathfrak{m}(x,\cdot)$ belong to the domain of the Friedrichs Laplacian and, hence, $\mathfrak{m}(x,y)=c(y)=\tilde{c}(x)={\rm const}$. One can assume (by specifying the solution to \eqref{correction robin eq}) that ${\rm const}=0$. Then formula \eqref{Green Verlinde} is valid. To prove that the above defined $m$ is indeed Robin mass \eqref{Robin mass}, it is sufficient to note that the near-diagonal asymptotics
\[\Phi(x,y)=-\frac{\log d(x,y)}{2\pi}+o(1), \qquad x\to y\in\dot{X},\]
is valid for function \eqref{Phi func}.
\end{proof}

\begin{Remark}\label{Ver Rem}
Clearly, $m\in C^{\infty}(\dot{X})$. From \eqref{Green Verlinde} and \eqref{Phi near vertices}, it follows that $m$ has logarithmic singularities near vertices
\[m(x)=\frac{b_k}{\beta_k}\log d(x,P_k)+O(1), \qquad x\to P_k.\]
\end{Remark}
\begin{Remark}\label{smoth conical case rem 1}
Formula~\eqref{Green Verlinde} holds if the conical metric $\mathfrak{T}$ on $X$ is flat in {\rm(}punctured{\rm)} vicinities of the conical points and is arbitrary (smooth) outside these vicinities. The proof of this fact is the same as in Lemma~\ref{Ver Vert Lem} except for the following replacements: 1) formula \eqref{delta of phi} is replaced by
\[\Delta_{x}\bigl(G(x,y)-\Phi(x,y)\bigr)=\mathfrak{T}^{-1}(x)B(x,x)-A^{-1}+K(x)/4\pi,\]
where $K=2\mathfrak{T}\partial\overline{\partial}\log \mathfrak{T}$ is the Gaussian curvature of $\mathfrak{T}$; 2) the function $\tilde{m}$ is defined as a solution to the equation
\[\Delta\tilde{m}(z)=\mathfrak{T}^{-1}B(z,z)-A^{-1}-\mathfrak{F}(z)+K(x)/4\pi\]
instead of {\rm \eqref{correction robin eq}}. The latter equation is solvable due to \eqref{for ortho} and the Gauss--Bonnet formula
\[\int_X K\,{\rm d }S-\sum_{k=1}^M b_k=2-2g.\]
\end{Remark}

\section[Generalization of Steiner formula (1.4) for polyhedra]{Generalization of Steiner formula (\ref{Zeta Robin}) for polyhedra}\label{Steiner sec}
In this section, we prove the following statement.
\begin{Proposition}
\label{Zeta_Robin_prop}
Formula {\rm \eqref{Zeta Robin}} is valid for a Troyanov metric $\mathfrak{T}$.
\end{Proposition}

As in \cite{Steiner}, the proof of Proposition \ref{Zeta_Robin_prop} is based on the short-time asymptotics of the heat kernel $H$ for the (Friedrichs) Laplacian $\Delta$. For this reason, we first derive (in Sections \ref{loc sol sec}--\ref{hk asym rem ssec}) such an asymptotics; it includes the new terms arising due the presence of conical singularities of the metric.

For general manifolds with conical singularities and Laplacians acting on forms, a similar asymptotics (including short-time asymptotics for the heat trace) is constructed in \cite{Mooers}.
 Although we adopt the scheme from \cite{Mooers}, in order to prove (4) we use a more explicit construction that provides the exponential decay of the remainder (including the weighted estimates near vertices). In particular, we fill a certain gap between the explicit constructions in dimension two based on Sommerfeld-type integrals (going back to~\cite{Carslaw} and often used by physicists, see, e.g.,~\cite{Dowker}) and the framework of the operator theory.

\subsection[Local solutions to the heat equation (partial\_t+Delta)u=0]{Local solutions to the heat equation $\boldsymbol{ (\partial_t+\Delta)u=0}$}
\label{loc sol sec}
\subsection*{Local solutions outside conical points} Let $U\subset X$ be a domain containing no conical points and let $z$ be a holomorphic coordinate on~$U$ in which $\mathfrak{T}=|{\rm d}z|^2$. Then $\Delta=-4\partial_z\partial_{\overline{z}}$ and the local fundamental solution to $(\partial_t+\Delta)u=0$ in $U$ is given by
\begin{equation}
\label{HK Plane}
\mathcal{H}(P,P',t|U,z)=\frac{1}{4\pi t}\exp\left(-\frac{|z-z'|^2}{4t}\right),
\end{equation}
where $z=z(P)$, $z'=z(P')$ and $P,P'\in U$.

\subsection*{Local solutions near conical point} Now let $U\subset X$ be a neighbourhood containing one conical point $Q$ (of conical angle $\beta=2\pi({b+1})$, $b>-1$) and let $z$ be a holomorphic coordinate on $U$ obeying $z(Q)=0$ and $\mathfrak{T}=\big|{\rm d}\bigl(z^{\beta/2\pi}\bigr)\big|^2$. Introduce the polar coordinates
\begin{equation}
\label{polar coordinates}
r=|z|^{\beta/2\pi}, \qquad \varphi=\frac{\beta}{2\pi}\operatorname{arg}z \quad ({\rm mod}\ \beta).
\end{equation}
Then $\Delta=-r^{-2}\bigl((r\partial_r)^2+\partial_{\varphi}^2\bigr)$ and the local fundamental solution to $(\partial_t+\Delta)u=0$ near $Q$ is given by
\begin{equation}
\label{HK Cone}
\mathcal{H}(P,P',t|U,z)=\frac{1}{8\pi {\rm i} \beta t}\int_{\mathscr{C}}\exp\left(-\frac{\mathfrak{r}^2}{4t}\right)\operatorname{cot}\Theta \, {\rm d}\vartheta
\end{equation}
(see \cite{Carslaw}, see also \cite{Dowker,KokKorJDG}), where
\begin{equation}\label{rtheta}
\mathfrak{r}:=\sqrt{r^2-2rr'\cos \vartheta+r'^{2}}, \qquad \Theta:=\pi\beta^{-1}(\vartheta+\varphi-\varphi')
\end{equation}
and $r=r(P)$, $\varphi=\varphi(P)$, $r'=r(P')$, $\varphi'=\varphi(P')$. The contour $\mathscr{C}$ is given by $\mathscr{C}=\mathscr{C}_+\cup\mathscr{C}_-$, where~$\mathscr{C}_\pm$ is a curve in the semi-strip $\{\operatorname{Re}\theta\in(-\pi,\pi),\pm\operatorname{Im}\theta>0\}$ going from $\pm ({\rm i}\infty+\pi)$ to~${\pm({\rm i}\infty-\pi)}$. The contour $\mathscr{C}$ in \eqref{HK Cone} can be replaced by the union, $\mathfrak{C}$, of the lines $\pm l:=\{\vartheta=\pm(\pi-{\rm i}\vartheta')\}_{\vartheta'\in\mathbb{R}}$ and small anti-clockwise circles $\odot[\vartheta_*]$ centred at the roots $\vartheta_*$ of $\operatorname{cot}\Theta$. Denote $\tilde{\mathfrak{C}}:=\mathfrak{C}\backslash\odot[\varphi'-\varphi]$, then \eqref{HK Cone} can be rewritten as $\mathcal{H}=\mathcal{H}_0+\tilde{\mathcal{H}}$, where
\begin{gather}
\mathcal{H}_0(P,P',t|U,z)=\frac{1}{4\pi t}\exp\left(-\frac{\mathfrak{d}^2}{4t}\right), \nonumber\\
\tilde{\mathcal{H}}(P,P',t|U,z)=\frac{1}{8\pi {\rm i} \beta t}\int_{\tilde{\mathfrak{C}}}\exp\left(-\frac{\mathfrak{r}^2}{4t}\right)\operatorname{cot}\Theta\, {\rm d}\vartheta,\label{HK Cone 1}
\end{gather}
and \smash{$\mathfrak{d}:=\bigl(r^2-2rr'\cos (\varphi'-\varphi)+r'^{2}\bigr)^{1/2}$} coincides with the distance
\[
d(P,P') \quad \text{if} \ d(P,P')\le \min\{d(P,Q),d(P',Q)\}.\]
 Since $\mathfrak{r}^2=(r-r')^{2}+4rr'\cosh ^2(\vartheta'/2)$ on $\pm l$, the integrand in \eqref{HK Cone 1} decays super-exponentially as~${|\vartheta|\to\infty}$. Therefore, $\tilde{\mathcal{H}}(P,P',t|U,z)$ decays exponentially as~${t\to 0}$ uniformly in $P$, $P'$ separated from $Q$.

\subsection{Parametrix construction}
\label{parametrix ssec}
Let $\{(U_k,z_k)\}_{k=1}^K$ be a complex atlas on $X$ and let $\{\chi_k\}_{k=1}^K$ be a partition of unity subordinate to this atlas. We assume that, for $k=1,\dots,M$, the relations $U_k\ni P_k$ and $\chi_k(P)=\chi(r_k(P))$ ($P\in U_k$) are valid, where $\chi$ is a smooth cut-off function equal to 1 near the origin and such that the diameter of $\operatorname{supp}\chi$ is sufficiently small, $(r_k,\varphi_k)$ are given by \eqref{polar coordinates} with $z=z_k$. Introduce the smooth cut-off functions $\{\kappa_k\}_{k=1}^K$ in such a way that $\operatorname{supp} \kappa_k\subset U_k$ and $\kappa_k=1$ in some neighborhood of $\operatorname{supp}\chi_k$.

Let us take
\begin{equation}
\label{parametrix}
\mathscr{H}^{(0)}(P,P',t)=\sum_{k=1}^K\kappa_k(P)\mathcal{H}(P,P',t|U_k,z_k)\chi_k(P')
\end{equation}
as a parametrix for the heat kernel $H$ as $t\to+0$. To show that $\mathscr{H}^{(0)}(\cdot,P',t)$ belongs to the domain of the Friedrichs Laplacian $\Delta$ for any $P'\in\dot{X}$, $t>0$, it is sufficient to show that the asymptotics of local solution \eqref{HK Cone} as $r\to 0$ does not contain growing terms. Indeed, using the expansion
\begin{align*}
\operatorname{cot}\Theta=\mp 2{\rm i}\left(\frac{1}{2}+\sum_{k=1}^\infty \exp\bigl(\pm 2k {\rm i}\Theta\bigr)\right), \qquad \vartheta\in\mathscr{C}_\pm,
\end{align*}
one can rewrite \eqref{HK Cone} as
\begin{gather}
\beta t \mathcal{H}(P,P',t|U,z) \exp\left(\frac{r^2+r'^{2}}{4t}\right)\nonumber\\
\qquad=\frac{1}{4\pi {\rm i} }\int_{-{\rm i}\mathscr{C}_+}\exp\left(\frac{rr'\cosh \vartheta}{2t}\right) {\rm d}\vartheta\nonumber
\\
\phantom{\qquad=}{}+\sum_{k=1}^\infty \cos \left(\frac{2\pi k(\varphi-\varphi')}{\beta}\right)\frac{1}{2\pi {\rm i} }\int_{-i\mathscr{C}_+}\exp\left(\frac{rr'\cosh \vartheta}{2t}-\frac{2\pi k\vartheta}{\beta}\right){\rm d}\vartheta\nonumber\\
\qquad=\frac{1}{2}I_0\left(\frac{rr'}{2t}\right)+\sum_{k=1}^\infty \cos \left(\frac{2\pi k(\varphi-\varphi')}{\beta}\right)I_{2\pi k/\beta}\left(\frac{rr'}{2t}\right),\label{boundedness of heat k near v}
\end{gather}
where
\[I_{\nu}(z)=\frac{1}{2\pi {\rm i} }\int_{-i\mathscr{C}_+}\exp(z\cosh \vartheta-\nu\vartheta)\, {\rm d}\vartheta= z^{\nu}\sum_{q=0}^{\infty}\iota_k(\nu)z^{2k}\]
is the modified Bessel function of the first kind. Taking into account the super-exponential decay of the coefficients $\iota_k(\nu)=1/q!2^{2q+\nu}\Gamma(q+\nu+1)$ of $I_{\nu}$ as $k,\nu\to +\infty$, one can majorize the absolute value of the right-hand side of~\eqref{boundedness of heat k near v} by the converging series $\sum_{k,q}c_k(rr'/2t)^{2\pi k/\beta+2q}=O(1)$
as~${r\to 0}$ with fixed $r', t>0$.
So, the \smash{$\mathscr{H}^{(0)}(\cdot,P',t)$} is bounded near the vertices and belongs to the domain, $\operatorname{Dom}\Delta$, of the Friedrichs Laplacian $\Delta$,
\[\operatorname{Dom}\Delta=\left\{u=\sum_{k=1}^M\chi_k(P)\sum_{l=-[b_k+1]}^{[b_k+1]}c_{k,l} r_k^{\frac{|l|}{b_k+1}}{\rm e}^{\frac{{\rm i}\varphi_k l}{b_k+1}}+\tilde{u} | c_{k,l}\in\mathbb{C},\, \tilde{u}\in\operatorname{Dom}\Delta_0\right\}\]
(see, e.g., \cite{Mooers}). Here $\operatorname{Dom}\Delta_0$ is the domain of the closure of~$\Delta_0$ in $L_2(X,\mathfrak{T})$, i.e., the closure of~$C_0^\infty(\dot{X})$ in the graph-norm \smash{$\|\tilde{u}\|^2_{\Delta_0}:=\|\tilde{u}\|^2_{L_2(X,\mathfrak{T})}+\|\Delta_0\tilde{u}\|^2_{L_2(X,\mathfrak{T})}$}.

Now substitution of \eqref{parametrix} into the heat equation yields
\begin{equation}
\label{discrepancy}
(\partial_t +\Delta)\mathscr{H}^{(0)}(\cdot,P',t)=\sum_{k=1}^K\chi_k(P') [\Delta,\kappa_k]\mathcal{H}_k(\cdot,P',t|U_k,z_k)=:\mathfrak{F}(\cdot,P',t)
\end{equation}
in $L_2(X,\mathfrak{T})$, where $t>0$ and $P'\in\dot{X}$. Since $[\Delta,\kappa_k]=0$ in a neighbourhood of $\operatorname{supp}\chi_k$, we have~${d(P,P')\ge {\rm const}>0}$ for any $(P,P',t)\in \operatorname{supp}\mathfrak{F}$. Due to \eqref{HK Plane} and \eqref{HK Cone}, this means that $\mathfrak{F}(\cdot,\cdot|t)$ and all its partial derivatives (with respect to the above local coordinates) decay exponentially as $t\to +0$.

Since $\mathcal{H}(P,P',t|U,z)\to\delta(z-z')$ as $t\to +0$ holds (in the sense of distributions) for solutions~\eqref{HK Plane}, \eqref{HK Cone}, we have
\begin{equation}
\label{delta convergence}
\mathscr{H}^{(0)}(\cdot,P',t)\to \sum_{k=1}^K(\kappa_k\chi_k)(P')\delta_{P'}=\sum_{k=1}^K\chi_k(P')\delta_{P'}=\delta_{P'}, \qquad t\to+0,
\end{equation}
where $\delta_{P'}$ is the Dirac measure at $P'$.

Consider the Cauchy problem
\begin{equation}
\label{remainder eq}
(\partial_t +\Delta)\tilde{\mathscr{H}}(\cdot,P',t)=-\mathfrak{F}(\cdot,P',t), \qquad t>0, \qquad \tilde{\mathscr{H}}(\cdot,P',0)=0.
\end{equation}
For any $P'\in X$, the solution $\tilde{\mathscr{H}}(\cdot,P',\cdot)\in C^{1}(\mathbb{R}_+;L_2(X,\mathfrak{T}))\cap C(\mathbb{R}_+;\operatorname{Dom}\Delta)$ to \eqref{remainder eq}
is given~by
\begin{gather}
\label{remainder hk}
\tilde{\mathscr{H}}(\cdot,P',t)=\int_0^t\exp\bigl((t'-t)\Delta\bigr)[-\mathfrak{F}(\cdot,P',t')]\,{\rm d}t'.
\end{gather}
Here $\tau\mapsto \exp(-\tau\Delta)$ ($\tau\ge 0$) is the $C_0$-semigroup with generator $\Delta$.

Denote
\begin{equation}
\label{heat kernel}
H=\mathscr{H}^{(0)}+\tilde{\mathscr{H}};
\end{equation}
then $H(\cdot,P',\cdot)\in C^{1}((0,+\infty);L_2(X,\mathfrak{T}))\cap C((0,+\infty);\operatorname{Dom}\Delta)$ and equations \eqref{discrepancy} and \eqref{remainder eq} imply $(\partial_t +\Delta)\mathscr{H}(\cdot,P',t)=0$ for $P'\in X$ and $t>0$.

Now we are to show that $H$ is the heat kernel of $\Delta$ (let us emphasize that the methods we use here are fairly standard for the smooth setting, see, e.g., \cite{BGV}).

 Let $\{\lambda_k\}_{k=0}^{\infty}$ be the sequence of all the eigenvalues of $\Delta$ counted with their multiplicities and let $\{u_k\}_{k=0}^{\infty}$ be the orthonormal basis of the corresponding eigenfunctions. Denote
\[c_k(t,P')=(\mathscr{H}(\cdot,P',t),u_k)_{L_2(X,\mathfrak{T})};\]
then
\[-\partial_t c_k=-(\partial_t\mathscr{H},u_k)_{L_2(X,\mathfrak{T})}=(\Delta\mathscr{H},u_k)_{L_2(X,\mathfrak{T})}=(\mathscr{H},\Delta u_k)_{L_2(X,\mathfrak{T})}=\lambda_k c_k.\]
Next, since $\Delta$ is non-negative, the operator $\exp((t'-t)\Delta)$ is a subcontraction for $t'<t$. Then formula \eqref{remainder hk} and estimates for $\mathfrak{F}$ given after \eqref{discrepancy} imply that the $L_2(X,\mathfrak{T})$-norm of $\tilde{\mathscr{H}}(\cdot,P',t)$ decays exponentially and uniformly in $P'\in X$ as $t\to +0$. Now, \eqref{delta convergence} implies
\begin{align*}
c_k(t,P')=&(H(\cdot,P',t),u_k)_{L_2(X,\mathfrak{T})}\\
=&\big(\mathscr{H}^{(0)}(\cdot,P',t),u_k\big)_{L_2(X,\mathfrak{T})}+o(1)\to \overline{u_k(P')}, \qquad t\to +0.
\end{align*}
Therefore, \smash{$c_k(t,P')={\rm e}^{-\lambda_k t}\overline{u_k(P')}$} and \smash{$H(P,P',t)=\sum_{k=0}^{\infty}{\rm e}^{-\lambda_k t}u_k(P)\overline{u_k(P')}$} is indeed the heat kernel of $\Delta$.

\subsection{Estimate of the remainder}
\label{hk asym rem ssec}
Let us estimate the remainder $\tilde{\mathscr{H}}$ in \eqref{heat kernel}. As shown above, the $L_2(X,\mathfrak{T})$-norm of $\tilde{\mathscr{H}}(\cdot,P',t)$ decays exponentially and uniformly in $P'\in X$ as $t\to +0$. Denote
\[\dot{\tilde{\mathscr{H}}}(\cdot,P',t)=\int_0^t\exp\bigl((t'-t)\Delta\bigr)[-\partial_{t'}\mathfrak{F}(\cdot,P',t')]\,{\rm d}t' .\]
Then the function
\[(x,P',t)\mapsto\int_0^t \dot{\tilde{\mathscr{H}}}(x,P',t')\, {\rm d}t'\] belongs to $C^{2}(\mathbb{R}_+;L_2(X,\mathfrak{T}))\cap C^{1}(\mathbb{R}_+;\operatorname{Dom}\Delta)$ and obeys \eqref{remainder eq}. Therefore, it coincides with~$\tilde{\mathscr{H}}$ due to the uniqueness of the solution to Cauchy problem \eqref{remainder eq}. In particular, we have \smash{$\dot{\tilde{\mathscr{H}}}=\partial_t \dot{\tilde{\mathscr{H}}}$} and
\smash{$-\Delta\tilde{\mathscr{H}}=\dot{\tilde{\mathscr{H}}}(\cdot,P',t)+\mathfrak{F}$}.
Since $\mathfrak{F}(\cdot,\cdot|t)$, $\partial_t\mathfrak{F}(\cdot,\cdot|t)$ decay exponentially as $t\to +0$ and uniformly in $P'\in X$ and $\exp((t'-t)\Delta)$ is a subcontraction for $t'<t$, the $L_2(X,\mathfrak{T})$-norm of~$\Delta\tilde{\mathscr{H}}(\cdot,P',t)$ decays exponentially and uniformly in $P'\in X$.

Next, we make use of the following statement.
\begin{Lemma}
Denote by $d_{\rm ver}(x)$ the geodesic distance between $x$ and $\{P_1,\dots,P_M\}$ in the metric~$\mathfrak{T}$. For $\epsilon>0$, the estimate
\begin{equation}
\label{inc smooth}
\|d_{\rm ver}^{\epsilon}u\|_{C(X)}\le c(X,\mathfrak{T},\epsilon)(\|\Delta u\|_{L_2(X,\mathfrak{T})}+\|u\|_{L_2(X,\mathfrak{T})}), \qquad u\in \operatorname{Dom}\Delta,
\end{equation}
holds. If all the coefficients $b_k$ of the divisor of conical points are non-integer then {\rm\eqref{inc smooth}} is also valid for $\epsilon=0$.
\end{Lemma}
\begin{proof}
Using the notation of Section \ref{parametrix ssec}, put \smash{$\tilde{U}=\bigcup_{k>M}U_k$} and \smash{$\tilde{\chi}\!=\!\sum_{k>M}\chi_k$}; then ${\operatorname{supp}\tilde{\chi}\!\subset\!\tilde{U}}$ and the closure of $\tilde{U}$ does not contain vertices. Since the metric $\mathfrak{T}$ is smooth on $\tilde{U}$, the estimate%
\begin{equation}\label{inc smooth outside vert}
\|u\|_{C(\tilde{U})}\le c(X,\mathfrak{T},\tilde{U})(\|\Delta u\|_{L_2(\tilde{U},\mathfrak{T}|_{\tilde{U}})}+\|u\|_{L_2(\tilde{U},\mathfrak{T}|_{\tilde{U}})}), \qquad u\in \operatorname{Dom}\Delta,
\end{equation}
follows from the elliptic regularity theorem and Sobolev embedding theorem.

Therefore, it is sufficient to prove \eqref{inc smooth} with $X$ replaced by a small neighbourhood ${U=U_k}$ of a single vertex $Q=P_k$ ($k=1,\dots,M$). For simplicity, we assume for a while that $\operatorname{supp}u\subset U$ and $\Delta u$ is smooth. To prove the claim, we use the methods of the theory of elliptic equations in domains with piece-wise smooth boundaries \cite{Kon,NP}. Rewriting the Laplacian in polar coordinates \eqref{polar coordinates}, we obtain
\[\Delta=-r^{-2}\bigl((r\partial_r)^2+\partial_\varphi^2\bigr)=-{\rm e}^{-2\sigma}\bigl(\partial_\sigma^2+\partial_\varphi^2\bigr),\]
where $\sigma=\log r$. Denote $f=-r^2\Delta u$, then
\begin{equation}
\label{ellip eq 1}
\bigl(\partial_\sigma^2+\partial_\varphi^2\bigr)u=f.
\end{equation}
Introduce the Kondrat'ev norms (see \cite{Kon})
\begin{align*}%\label{weighted norm}
\|v\|_{H_\nu^l}&=\left(\sup_{\phi_0}\sum_{p+q\le l}\int_{0}^{+\infty}r{\rm d}r\int_{\phi_0}^{\phi_0+\beta} {\rm d}\varphi |\partial^q_\varphi(r\partial_r)^p v(r,\varphi)|^2 r^{2(\nu-1)}\right)^{1/2}\\
&=\left(\sup_{\phi_0}\sum_{p+q\le l}\int_{-\infty}^{+\infty}{\rm d}s \int_{\phi_0}^{\phi_0+\beta} {\rm d}\varphi |\partial^q_\varphi\partial_s^p v(s,\varphi)|^2 {\rm e}^{2\nu s}\right)^{1/2},
\end{align*}
where $\beta$ is the conical angle at $Q$, and $l=0,1,\dots$, $\nu\in\mathbb{R}$. The inclusions $\Delta u,u\in L_2(U,\mathfrak{T})$ imply $f\in H^0_{-1}$, $u\in H^0_1$. Moreover, since $\Delta u$ is smooth, we have $f\in H^0_{\nu'}$ for any $\nu'>-2$.

Let $\chi,\psi\in C_c^\infty(\mathbb{R})$ be smooth cut-off functions obeying $\psi(s)=1$ for~${s\in \operatorname{supp}\chi}$, $\chi(s)=1$ for $[-1,1]$, and $\psi(s)=0$ for $s\not\in[-2,2]$. The estimate
\[\|\chi u\|_{H_\nu^2}\le c(\|\psi u\|_{H_\nu^0}+\|\psi f\|_{H_\nu^0})\]
follows from the elliptic regularity theorem. Denote $\chi_L(s):=\chi(s-L)$ and $\psi_L(s):=\psi(s-L)$. Substituting $u_L(s)=u(s+L)$, $f_L(s)=f(s+L)$ in the above estimate and passing to the new variable $\sigma=s+L$, one obtains
\[\|\chi_L u\|_{H_\nu^2}\le c(\|\psi_L v\|_{H_\nu^0}+\|\psi_L f\|_{H_\nu^0}).\]
Summation of the above estimates with integer $L$ yields
\smash{$\|u\|_{H_\nu^2}\le c(\|u\|_{H_\nu^0}+\|f\|_{H_\nu^0})$}.
As a~corollary, the inclusions $f\in H^0_{-1}\subset H^0_1$, $u\in H^0_1$ imply $u\in H^2_\nu$ for any $\nu\ge 1$.

Applying the complex Fourier transform
\[\hat{v}(\mu,\varphi):=\frac{1}{\sqrt{2\pi}}\int_{-\infty}^{+\infty}{\rm e}^{-{\rm i}\mu s}v(s,\varphi)\,{\rm d}s,\]
one rewrites \eqref{ellip eq 1} as
\begin{equation}
\label{ellip eq 2}
\bigl(\partial_\varphi^2-\mu^2\bigr)\hat{u}(\varphi,\mu)=\hat{f}(\varphi,\mu), \qquad \varphi\in\mathbb{R}/\beta\mathbb{Z},
\end{equation}
where $\operatorname{Im}\mu=\nu$. Note that $\hat{f}$ is well-defined and holomorphic in $\mu$ for $\operatorname{Im}\mu>-2$ due to the inclusion $f\in H^0_{\nu'}$ ($\nu'>-2$). In view of the Parseval identity, we have the following equivalence of the norms
\[%\label{Parseval equiv}
\|v\|_{H_\nu^l}\asymp\left(\sup_{\phi_0}\sum_{p+q\le l}\int_{\nu i-\infty}^{\nu i+\infty}|\mu|^{2p}\|\partial^q_\varphi\hat{v}(\mu,\varphi)\|_{L_2([\phi_0,\phi_0+\beta])}^2 \,{\rm d}\mu\right)^{1/2}.
\]
The resolvent kernel for the equation \eqref{ellip eq 2} is given by
\[\mathfrak{R}_\mu(\varphi,\varphi')=-\frac{\cosh \bigl(\mu(|\varphi-\varphi'|-\beta/2)\bigr)}{2\mu \sinh(\beta\mu/2)}, \qquad \varphi,\varphi'\in [0,\beta].\]
The poles of the resolvent are $\mu=\frac{2\pi {\rm i} k}{\beta}$ with integer $k$. The estimate
\begin{equation}
\label{RK estimatesat infinity}
|\mu|^2\|\mathfrak{R}_\mu(\cdot,\cdot)\|_{L_2([\phi_0,\phi_0+\beta]^2)}+|\mu|\|\partial_\varphi\mathfrak{R}_\mu(\cdot,\cdot)\|_{L_2([\phi_0,\phi_0+\beta]^2)}=O(1)
\end{equation}
is valid for large $|\operatorname{Re}\mu|$.

The solution to \eqref{ellip eq 2} is given by
\[\hat{u}(\varphi,\mu)=\int_{0}^{\beta}\mathfrak{R}_\mu(\varphi,\varphi')\hat{f}(\varphi',\mu)\,{\rm d}\varphi', \]
Consider two solutions $u\in H^2_{\nu}$ and $v\in H^2_{\nu'}$ to \eqref{ellip eq 1}, where $-2<\nu'<\nu$, $\nu\ge 1$, and $\nu$, $\nu'$ do not coincide with any $\frac{2\pi k}{\beta}$ with integer $k$. Then the residue theorem implies
\begin{align}
u(r,\varphi)=&\frac{1}{2\pi}\int_{\nu {\rm i}-\infty}^{\nu i+\infty}{\rm e}^{{\rm i}\mu s}\int_{0}^{\beta}\mathfrak{R}_\mu(\varphi,\varphi')\int_{-\infty}^{+\infty}{\rm e}^{-{\rm i}\mu s'}f(s',\varphi')\, {\rm d}s' {\rm d}\varphi'{\rm d}\mu\nonumber\\
=&-{\rm i}\sum_{k}\int_{-\infty}^{+\infty}{\rm d}s'\int_{0}^{\beta}{\rm d}\varphi'f(s',\varphi')\underset{\mu=-2\pi {\rm i} k/\beta}{\rm Res}\bigl(\mathfrak{R}_\mu(\varphi,\varphi'){\rm e}^{{\rm i}\mu (s-s')}\bigr)+v(r,\varphi)\nonumber\\
=&-\sum_{k}\sum_{\pm}\mathcal{Y}_{k,\pm}(r,\varphi)(\Delta u,\mathcal{Y}_{k,\pm})_{L_2(U,\mathfrak{T}|_U)}+v(r,\varphi),\label{expansion near vertex}
\end{align}
where the summation is over all $k\in (-\nu\beta/2\pi,-\nu'\beta/2\pi)$ and
\begin{gather*}
\mathcal{Y}_{k,\pm}(r,\varphi)=\frac{1}{\sqrt{4\pi |k|}}r^{\frac{2\pi}{\beta}k}{\rm e}^{\pm {\rm i}\frac{2\pi}{\beta}k\varphi} ,\qquad k\ne 0,\\
\mathcal{Y}_{0,+}(r,\varphi)=\frac{1}{\sqrt{\beta}}, \qquad \mathcal{Y}_{0,-}(r,\varphi)=\frac{1}{\sqrt{\beta}}\log r.
\end{gather*}

Put some $\nu'\in (-1,-2)$ in \eqref{expansion near vertex}. Then the inclusion $v\in H^2_{\nu'}$ implies $v|_U\in H^2(U)\subset C(U)$. At the same time, $u$ is bounded near the vertex due to the inclusion $u\in \operatorname{Dom}\Delta$ and the assumption that $\Delta u$ is smooth. Therefore, asymptotics \eqref{expansion near vertex} cannot contain the (growing near the vertex) terms $\mathcal{Y}_{k,\pm}$, $k<0$ or $\mathcal{Y}_{0,-}$, i.e., the coefficients $(\Delta u,\mathcal{Y}_{k,\pm})_{L_2(U,\mathfrak{T}|_U)}$, $k<0$ and~${(\Delta u,\mathcal{Y}_{0,-})_{L_2(U,\mathfrak{T}|_U)}}$ are equal to zero.

Now, put $\nu'=-1+\epsilon$ into \eqref{expansion near vertex}, where $\epsilon>0$ can be arbitrarily small (or equal to zero if~$\beta/2\pi$ is non-integer). Then the inclusion $v\in H^2_{\nu'}$ implies $r^\epsilon v\in H^2(U)$ and
\[\|r^\epsilon v\|_{C(U)}\le c\|r^\epsilon v\|_{H^2(U)}\le c\|v\|_{H^2_{\nu'}}\]
due to the Sobolev embedding theorem. In addition, estimate \eqref{RK estimatesat infinity} implies
\[\|v\|_{H^2_{\nu'}}\le c\|f\|_{H^0_{\nu'}}\le c\|\Delta u\|_{L_2(U,\mathfrak{T}|_U)}.\]
Since $\Delta u\mapsto (\Delta u,\mathcal{Y}_{k,\pm})_{L_2(U,\mathfrak{T}|_U)}$ (with $k>-\beta/2\pi$) are bounded linear functionals of $\Delta u\in L_2(U,\mathfrak{T}|_U)$, the two last estimates lead to
\begin{equation}
\label{inc smothness near vertex}
\|r^\epsilon u\|_{C(U)}\le c\|\Delta u\|_{L_2(U,\mathfrak{T}|_U)}, \qquad u\in \operatorname{Dom}\Delta, \qquad \operatorname{supp} u\subset U.
\end{equation}

Now, consider the general case $v\in \operatorname{Dom}\Delta$ and $\Delta v\in C^{\infty}(X)$. Then estimate \eqref{inc smooth outside vert} is valid for $u=\tilde{\chi}v$ and estimates \eqref{inc smothness near vertex} are valid for $u=\chi_k v$, $U=U_k$. Summing these estimates and applying the elliptic regularity theorem yields
\begin{align*}
\|d_{\rm ver}^{\epsilon}v\|_{C(X)}&\le c\left(\sum_{k=1}^M \|r^\epsilon_k \chi_k v\|_{C(U_k)}+\|\tilde{\chi}u\|_{C(\tilde{U})}\right)\\
&\le c\left(\|\Delta v\|_{L_2(X,\mathfrak{T})}
+\sum_{k=1}^M \|[\Delta,\chi_k ]v\|_{L_2(X,\mathfrak{T})}+\|[\Delta,\tilde{\chi}]v\|_{L_2(X,\mathfrak{T})}+\|v\|_{L_2(X,\mathfrak{T})}\right) \\
&\le c(\|\Delta v\|_{L_2(X,\mathfrak{T})}+\|v\|_{L_2(X,\mathfrak{T})}).
\end{align*}
Thus, inequality \eqref{inc smooth} is proved under the additional condition $\Delta u\in C^{\infty}(X)$. In the general case $\Delta v\in L_2(X,\mathfrak{T})$, there is a sequence $\{u_j\}$ converging to $u$ in the graph norm of $\Delta$ and such that $\Delta u_k\in C^{\infty}(X)$. Applying \eqref{inc smooth} to $u_j-u_{j'}$ we prove that $u_j\to u$ in the norm $\|d_{\rm ver}^{\epsilon}\cdot\|_{C(X)}$. It remains to take the limit as $j\to \infty$ in \eqref{inc smooth} with $u=u_j$.
\end{proof}

As a corollary of \eqref{inc smooth} and estimates made at the beginning of the subsection, for arbitrarily small $\epsilon>0$, the function $t\mapsto d_{\rm ver}^{\epsilon}(P)\tilde{\mathscr{H}}(P,P',t)$ decays exponentially and uniformly in $(P,P')\in X\times X$ as $t\to +0$.

Introduce the shorter notation $d=d(P,P')$. In what follows, we assume that $P$ and $P'$ are closer to each other than to the vertices, i.e., $d\le d_{\rm ver}(P)$. Now combining formulas \eqref{HK Plane}, \eqref{HK Cone 1}, \eqref{parametrix} and \eqref{heat kernel} leads to
\begin{equation}
\label{heat kernel expansion}
\begin{split}
\mathscr{H}(P,P',t)=\frac{1}{4\pi t}\exp\left(-\frac{d^2}{4t}\right)+\sum_{k=1}^M\tilde{\mathcal{H}}(P,P',t|U_k,z_k)\chi_k(P')+\tilde{\mathscr{H}}(P,P',t).
\end{split}
\end{equation}
Note that the sum in \eqref{heat kernel expansion} decays as $t\to +0$ exponentially and uniformly in $P$, $P'$ separated from the vertices.

\subsection{Formula for Robin mass}
For $\operatorname{Re} s>0$, the Green function $\mathscr{G}(\cdot,\cdot|s)$ for the operator $\Delta^{s}$ is related to the heat kernel via
\begin{equation}
\label{green func}
\mathscr{G}(P,P'|s)=\frac{1}{\Gamma(s)}\int_{0}^{+\infty}\bigl(\mathscr{H}(P,P',t)-A^{-1}\bigr)t^{s-1}\,{\rm d}t.
\end{equation}
Substitution of \eqref{heat kernel expansion} and \eqref{HK Cone 1} into \eqref{green func} leads to
\begin{align}
\mathscr{G}(P,P',s)={}&\frac{1}{4\pi\Gamma(s)}\int_{0}^{1}\exp\left(-\frac{d^2}{4t}\right)t^{s-2}\, {\rm d}t +\sum_{k=1}^M\chi_k(P')\mathcal{E}(P,P',s|U_k,z_k)
+\tilde{\mathscr{G}}(P,P',s)\nonumber\\
={}&\frac{1}{4\pi\Gamma(s)}\left(\frac{d^2}{4}\right)^{s-1}\Gamma\left(1-s\Big|\frac{d^2}{4}\right)+\sum_{k=1}^M\chi_k(P')\mathcal{E}(P,P',s|U_k,z_k)\nonumber\\
&
+\tilde{\mathscr{G}}(P,P',s),\label{Green exp}
\end{align}
where
\[\Gamma(z|\tau)=\int_{\tau}^{+\infty}\exp(-u)u^{z-1}\, {\rm d}u\]
is the incomplete Gamma function and
\begin{gather}
\mathcal{E}(P,P',s|U,z)=\frac{1}{8\pi {\rm i} \beta\Gamma(s)}\int_{0}^{1} {\rm d}t\, t^{s-2}\int_{\tilde{\mathfrak{C}}}{\rm d}\vartheta \exp\left(-\frac{\mathfrak{r}^2}{4t}\right)\operatorname{cot}\Theta,\label{molel tilde green}
\\
\tilde{\mathscr{G}}(P,P'|s)=\frac{1}{\Gamma(s)}\left[\int_{1}^{+\infty}\bigl(\mathscr{H}(P,P',t)-A^{-1}\bigr)t^{s-1}\, {\rm d}t\right.\nonumber\\
\phantom{\tilde{\mathscr{G}}(P,P'|s)=}{}\left.+\int_{0}^{1}\tilde{\mathscr{H}}(P,P',t)t^{s-1}\, {\rm d}t-(sA)^{-1}\right].\label{tilde green}
\end{gather}

Since $\mathscr{H}(P,P',t)-A^{-1}$ and $\tilde{\mathscr{H}}(P,P',t)$ decay exponentially and uniformly in $P,P'\in X$ as~${t\to +\infty}$ and $t\to +0$, respectively, the function $\tilde{\mathscr{G}}(\cdot,\cdot|s)$ is smooth in $X\times X$ and is holomorphic in $s$ in the half-plane $\operatorname{Re} s>0$.

Similarly, $\mathcal{E}(\cdot,\cdot,s|U,z)$ is smooth for $r,r'>0$ and is holomorphic in $s$ for $\operatorname{Re} s>0$. Therefore, if $P$, $P'$ are separated from the vertices, then $\mathscr{G}(P,P'|s)$ is holomorphic in $s$ in a neighbourhood of $s=1$.

Due to the asymptotics
\begin{align*}
-\lim_{s\to 1, \operatorname{Re} s<1}\Gamma(1-s|\tau)&=\lim_{s\to 1, \operatorname{Re} s<1}\left[\int_{0}^{\tau}\frac{{\rm e}^{-u}-1}{u^{s}}\, {\rm d}u+\frac{\tau^{1-s}}{1-s}-\Gamma(1-s)\right]\\
&=\gamma+\log \tau+O(\tau)
\end{align*}
as $\tau\to +0$, expansion \eqref{Green exp} yields the following formula for the Robin mass
\begin{align}
m(P):={}&\lim_{P'\to P}\left[\mathscr{G}(P,P',|1)+\frac{\log d}{2\pi}\right]\nonumber\\
={}&\frac{2\log 2-\gamma}{4\pi}
+\sum_{k=1}^M\chi_k(P)\mathcal{E}(P,P,1|U_k,z_k)+\tilde{\mathscr{G}}(P,P,1),\label{Robin pointwise}
\end{align}
where $P\in\dot{X}$.

Recall that the Robin mass $m$ is smooth on $\dot{X}$ and has logarithmic (and, therefore, integrable) singularities at $P_1,\dots,P_M$ (see Remark \ref{Ver Rem}). Integration of both parts of \eqref{Robin pointwise} over $X$ leads~to
\begin{gather}\label{Robin int 1}
\int_{X}m\, {\rm d }S=\frac{2\log 2-\gamma}{4\pi}A+\sum_{k=1}^M\mathfrak{P}(U_k,z_k,\chi_k)+\int_{X}\tilde{\mathscr{G}}(P,P,1)\, {\rm d }S(P).
\end{gather}
where
\[\mathfrak{P}(U,z,\chi)=\int_0^{\infty}{\rm d}r\, r\chi(r)\int_0^{\beta}{\rm d}\varphi\, \mathcal{E}(P,P,1|U,z).\]

Introduce the function
\begin{gather}\label{countor contribution}
\mathfrak{I}(r,t|\beta)=\frac{1}{16\pi {\rm i} }\int_{\tilde{\mathfrak{C}}}\exp\left(-\frac{r^2\sin^2(\vartheta/2)}{t}\right)\frac{\operatorname{cot}(\pi\vartheta/\beta)}{{\rm sin}^2(\vartheta/2)}\, {\rm d}\vartheta
\end{gather}
and denote $\mathfrak{I}'=\partial_r\mathfrak{I}$ and $\dot{\mathfrak{I}}=\partial_t\mathfrak{I}$. The integrand in \eqref{countor contribution} decays super-exponentially as $|\vartheta|\to \infty$ (for any $t\ge 0$) and exponentially as $t\to +0$ and $r>0$. Then $\mathfrak{I}(r,t|\beta)$ is differentiable in $r$, $t$ for~${r>0}$ and is continuous for $r,t\ge 0$, $r+t>0$. Note that
\[\mathfrak{I}(0,t|\beta)=\frac{1}{12}\left(\frac{2\pi}{\beta}-\frac{\beta}{2\pi}\right), \qquad \mathfrak{I}(r,0|\beta)=0.\]
Formulas \eqref{molel tilde green}, \eqref{rtheta} and \eqref{countor contribution} lead to
\begin{align*}
\int_{R_-}^{R_+}{\rm d}r\, r\int_{0}^{\beta}{\rm d}\varphi\, \mathcal{E}(P,P,1|U,z)&=\int_{R_-}^{R_+}\frac{{\rm d}r r}{8\pi {\rm i} }\int_{0}^{1}\frac{{\rm d}t}{t}\int_{\tilde{\mathfrak{C}}}{\rm d}\vartheta \exp\left(-\frac{r^2 \sin^2(\vartheta/2)}{t}\right)\operatorname{cot}\left(\frac{\pi\vartheta}{\beta}\right)\\
&=-\int_{R_-}^{R_+}{\rm d}r\int_{0}^{1}{\rm d}t\, \mathfrak{I}'(r,t|\beta).
\end{align*}
Therefore,
\begin{gather}
\label{J alpha}
\mathfrak{P}(U,z,\chi)=-\int_{0}^{1}{\rm d}t\int_0^{\infty}{\rm d}r\, \chi(r)\mathfrak{I}'(r,t|\beta).
\end{gather}
Now we are ready to prove Proposition~\ref{Zeta_Robin_prop}.

\begin{proof}[Proof of Proposition \ref{Zeta_Robin_prop}] Recall that the $\zeta$-function is related to the trace of the heat kernel
\[\mathscr{H}(t):=\sum_{k=0}^{\infty}{\rm e}^{-\lambda_k t}=\int_X \mathscr{H}(P,P,t)\,{\rm d }S(P)\]
via
\[\zeta_\Delta(s)=\frac{1}{\Gamma(s)}\int_{0}^{\infty}\bigl(\mathscr{H}(t)-1\bigr)t^{s-1}\,{\rm d}t.\]
In view of \eqref{heat kernel expansion}, \eqref{HK Cone 1}, \eqref{rtheta}, and \eqref{countor contribution}, we have
\begin{gather}
\mathscr{H}(t)-\frac{A}{4\pi t}-\int_X\tilde{\mathscr{H}}(P,P,t)\,{\rm d }S(P)\nonumber\\
\qquad=\sum_{k=1}^M\int_X\tilde{\mathcal{H}}(P,P,t|U_k,z_k)\chi_k(P)\,{\rm d }S(P)\nonumber\\
\qquad=\frac{1}{8\pi {\rm i} t}\sum_{k=1}^M\int_{r=0}^{+\infty}{\rm d}r_k\, r_k\chi_k(r_k)\int_{\tilde{\mathfrak{C}}}{\rm d}\vartheta \exp\left(-\frac{r_k^2\sin^2(\vartheta/2)}{t}\right)\operatorname{cot}\left(\frac{\pi\vartheta}{\beta_k}\right)\nonumber\\
\qquad=-\sum_{k=1}^M\int_{0}^{\infty}\chi_k(r_k)\mathfrak{I}'(t,r_k,\beta_k)\,{\rm d}r_k.\label{spec part}
\end{gather}
Note that the right-hand side of \eqref{spec part} is continuous in $t\ge 0$. Formulas \eqref{reg zeta of 1}, \eqref{spec part}, and \eqref{countor contribution} and the lead to
\begin{gather*}
\operatorname{reg}\zeta_\Delta(1)=\lim_{s\to 1, \operatorname{Re} s>1}\left(\frac{1}{\Gamma(s)}\int_0^\infty (\mathscr{H}(t)-1 )t^{s-1}\, {\rm d}t-\frac{A}{4\pi(s-1)}\right)
\\
\qquad=\int_1^\infty (\mathscr{H}(t)-1 )\,{\rm d}t-1+\lim_{s\to 1, \operatorname{Re} s>1}\left(\frac{1}{\Gamma(s)}\int_0^1\mathscr{H}(t)t^{s-1}\, {\rm d}t-\int_0^1\frac{A}{4\pi t}t^{s-1}\,{\rm d}t\right)
\\
\qquad=\int_1^\infty (\mathscr{H}(t)-1 )\,{\rm d}t-1-\frac{A}{4\pi}\lim_{s\to 1, \operatorname{Re} s>1}\left(\left(\frac{1}{\Gamma(s)}-1\right)\int_0^1t^{s-2}\, {\rm d}t\right)
\\
\phantom{\qquad=}{}-\lim_{s\to 1}\left(\frac{1}{\Gamma(s)}\int_0^1\sum_{k=1}^M\int_{0}^{\infty}\chi_k(r_k)\mathfrak{I}'(t,r_k,\beta_k)\,{\rm d}r_k t^{s-1}{\rm d}t\right)
\\
\qquad=\int_{X}{\rm d }S(P)\int_1^\infty {\rm d}t\big(\mathscr{H}(P,P,t)-A^{-1}\big)-1+\int_X\, {\rm d }S(P)\int_0^1 {\rm d}t \, \tilde{\mathscr{H}}(P,P,t)
\\
\phantom{\qquad=}{}-\sum_{k=1}^M\int_0^{1}\, {\rm d}t\int_{0}^{\infty}{\rm d}r_k\, \chi_k(r_k)\mathfrak{I}'(t,r_k,\beta_k)+\frac{A\gamma}{4\pi}.
\end{gather*}
In view of \eqref{tilde green}, \eqref{J alpha} and \eqref{Robin int 1}, the last equality implies \eqref{Zeta Robin}:
\begin{align*}
\operatorname{reg}\zeta_\Delta(1)&=\int_X\tilde{\mathscr{G}}(P,P|1)\, {\rm d }S(P)+\frac{A\gamma}{4\pi}+\sum_{k=1}^M\mathfrak{P}(U_k,z_k,\chi_k)\\
&=\int_{X}m\, {\rm d }S-\frac{\log 2-\gamma}{2\pi}A.
\end{align*}
Proposition \ref{Zeta_Robin_prop} is proved.
\end{proof}

\begin{Remark}\label{smoth conical case rem 2}
Combining the above arguments with Steiner's proof of \eqref{Zeta Robin} for the case of smooth metric (see~\cite{Steiner}), one concludes that \eqref{Zeta Robin} is valid for metrics smooth outside conical singularities and flat near them. In this case, parametrix \eqref{parametrix} for~$\mathscr{H}$ should be slightly modified in the domain $\tilde{U}=\bigcup_{k>M}U_k$ (see \cite[formula~(51)]{Steiner}).
\end{Remark}

\section[Asymptotics of reg zeta\_Delta(1) for polyhedron with small pinch]{Asymptotics of $\boldsymbol{\operatorname{reg} \zeta_\Delta(1)}$ for polyhedron with small pinch}
\label{asymp sec}
In this section, we introduce the family $\varepsilon\mapsto X(\varepsilon)$ of degenerating polyhedral surfaces $X(\varepsilon)$ and prove asymptotic formulas \eqref{zeta 1 asymp}--\eqref{zeta_1 asymp_1} for them.

\subsection*{Family of degenerating surfaces} Let $X_\pm$ be polyhedral surfaces of genera $g_\pm$ endowed with Troyanov metrics $\mathfrak{T}_\pm$, respectively. Suppose that $P_\pm\in X_\pm$ do not coincide with the conical points of $X_\pm$ and chose holomorphic local coordinates $z_\pm$ near $P_\pm$ obeying $z_\pm(P_\pm)=0$ and $\mathfrak{T}_\pm(z_\pm)=1$.

Let $\varepsilon>0$ be a small parameter. Cut the surfaces $X_\pm$ along the segments
\[\mathscr{C}_\pm(\varepsilon):=\{Q_\pm\in X_\pm \mid z_\pm(Q_\pm)\in [0,\varepsilon]\}\]
and then make the cross-gluing of $X'_\pm:=X_\pm\backslash\mathscr{C}_\pm(\varepsilon)$ along the boundaries in a usual way. As a~result, we obtain the new polyhedral surface $X=X(\varepsilon)$ of genus $g=g_++g_-$ with the metric~${\mathfrak{T}=\mathfrak{T}(\varepsilon)}$ whose singularities are the conical points of $X_\pm$ and the two conical points $Q$ and $Q'$ (of conical angles equal to $4\pi$) obtained by identifying the points of $X'_\pm$ with coordinates~${z_\pm=0}$ and $z_\pm=\varepsilon$, respectively. The area $A$ of $(X,\mathfrak{T})$ is equal to $A_++A_-$, where $A_{\pm}$ is the area of~${(X_\pm,\mathfrak{T}_\pm)}$.

In what follows, we omit the dependence on $\varepsilon$ in the notation for objects related to the surface~${X=X(\varepsilon)}$. At the same time the objects related to $X_\pm$ (such as Laplacians, Green functions, Robin masses, etc) are marked by the lower indices $\pm$. Also, we consider $X'_\pm$ as domains in $X$.

Near $Q$, we introduce the new holomorphic coordinates $\xi$ and $\zeta$,
\begin{equation}
\label{coordinates near pinch}
\xi=\pm\sqrt{\frac{z_\pm-\varepsilon}{z_\pm}} \longleftrightarrow z_\pm=\frac{\varepsilon}{1-\xi^2}, \qquad \zeta=\frac{\xi-1}{\xi+1}.
\end{equation}
The map $P\mapsto \zeta(P)$ defines the biholomorphism of some neighbourhood of $Q$ (the ``pinching zone'') onto the annulus $|\zeta|\in \bigl(c_0\varepsilon,c_1\varepsilon^{-1}\bigr)$ in the Riemann sphere $\overline{\mathbb{C}}$ obeying $\zeta(Q)=1$ and~${\zeta(Q')=-1}$. At the same time, the map $P\mapsto \xi(P)$ is the biholomorphism of the pinching zone on the domain $|\xi\pm 1|\ge c\varepsilon$ in $\overline{\mathbb{C}}$ while $\xi(Q)=\infty$ and $\xi(Q')=0$. Along with coordinates $z_\pm$, $\zeta$, $\xi$, we will use arbitrary local holomorphic coordinates on $X_\pm$ outside small neighbourhoods of $P_\pm$; all the coordinates listed above will be called {\it admissible}.

\subsection*{Formal asymptotics of $\boldsymbol{\operatorname{reg}\zeta_\Delta(1)}$ as $\boldsymbol{\varepsilon\to 0}$}
 We start this subsection with a formal derivation of asymptotic formula \eqref{zeta 1 asymp} and then give its justification.

Choose a canonical basis of cycles $\{a_{k,\pm},b_{\pm,k}\}_{k=1,\dots,g_{\pm}}$ on $X_\pm$ in such a way that none of the cycles intersects a small neighbourhood of $P_\pm$. Let $\vec{v}_\pm=\{v_{k,\pm}\}_{k=1,\dots,g_\pm}$ be the basis of the normalized Abelian differentials on $X_\pm$ and let $\mathbb{B}_\pm$ be the matrix of $b$-periods of $v_{k,\pm}$. Then
\[\big\{a_{1,+},\dots,a_{g_+,+},a_{1,-},\dots,a_{g_-,-},b_{1,+},\dots,b_{g_+,+},b_{1,-},\dots,b_{g_-,-}\big\}\]
is a canonical basis of cycles on $X$. Let $\vec{\omega}:=\{w_{1}^{+},\dots,w_{g_+}^{+},w_{1}^{-},\dots,w_{g_-}^{-}\}$ be the corresponding basis of the normalized Abelian differentials on $X$ and let $\mathbb{B}$ be its matrix of $b$-periods. As shown in \cite[Lemmas 1 and 2]{KokTAMS}, the asymptotics
\begin{gather}
\label{abel diff asymp}
w_{k}^{\pm}(x)\sim\begin{cases}
v_{k,\pm}(x), & x\in X_\pm,\\
0, & x\in X_\mp,
\end{cases}
\end{gather}
hold as $\varepsilon\to 0$ outside small neighbourhood of $Q$. As a corollary, we have
\begin{equation}
\label{priod matr asymp}
\mathbb{B}\sim\left(\begin{matrix}
\mathbb{B}_+ & 0\\
0 & \mathbb{B}_-
\end{matrix}\right).
\end{equation}
Formal integration of \eqref{abel diff asymp} yields
\begin{gather}
\label{abel map asymp}
\mathscr{A}(x-y)=\int_y^x\vec{\omega}\sim\begin{cases}
(\mathscr{A}_+(x-y),0)^{\mathsf T}, & x,y\in X_+,\\
(0,\mathscr{A}_-(x-y))^{\mathsf T}, & x,y\in X_-,\\
(\mathscr{A}_+(x-P_+),\mathscr{A}_-(P_--y))^{\mathsf T}, & x\in X_+, y\in X_-,\\
(\mathscr{A}_+(P_+-y),\mathscr{A}_-(x-P_-))^{\mathsf T}, & x\in X_-, y\in X_+,
\end{cases}
\end{gather}
where $\mathscr{A}$ and $\mathscr{A}_\pm$ are Abel maps for surfaces $X$ and $X_\pm$, respectively.

The asymptotics (see \cite[formulas (2.33) and (2.34)]{KokTAMS})
\begin{align}
\label{prime form asymp}
E(x,y)^2\sim\begin{cases}
E_\pm(x,y)^2, & x,y\in X_\pm, \\
-16\varepsilon^{-2}E_\pm^2(x,z_\pm=0)E^2_\mp(z_\mp=0,y), & x\in X_\pm, \ y\in X_\mp,
\end{cases}
\end{align}
holds (outside a small neighbourhood of $Q$) for the prime form $E$ of $X$; here $E_\pm$ are the prime forms of $X_\pm$.

The substitution of \eqref{abel map asymp} and \eqref{prime form asymp} into \eqref{Verlinde} and \eqref{Phi func} leads to
\begin{equation}
\label{Phi func asymp}
\Phi(x,y)\sim\begin{cases}
\Phi_\pm(x,y), & x,y\in X_\pm,\\
\Phi_\pm(x,P_\pm)+\Phi_\mp(y,P_\mp)+\frac{1}{2\pi}\log (\varepsilon/4), & x\in X_\pm, y\in X_\mp.
\end{cases}
\end{equation}
Formal integration of \eqref{Phi func asymp} over $X\times X$ and taking into account \eqref{scalar Robin} yield
\begin{align}
A\int_X m\, {\rm d }S&=-\sum_{\pm}\left[\int_{X_\pm\times X_\pm}+\int_{X_\pm\times X_\mp}\right]\Phi(x,y)\,{\rm d }S(x){\rm d }S(y)\nonumber\\
&\sim\sum_{\pm}\left(A_\pm\int_{X_\pm}m_\pm \,{\rm d }S-2A_\mp\int_{X_\pm}\Phi_\pm(\cdot,P_\pm)\,{\rm d }S\right)-\frac{A_+A_-}{\pi}\log (\varepsilon/4).\label{illegal int}
\end{align}
Integration of \eqref{Green Verlinde} with respect to $y$ leads to
\[-2\int_{X_\pm}\Phi(\cdot,P_\pm)\,{\rm d }S=A_\pm m_\pm(P_\pm)+\int_{X_\pm}m_\pm \, {\rm d }S.\]
Substitution of the last formula with $X$ replaced by $X_\pm$ into \eqref{illegal int} yields
\begin{equation}
\label{int Robin mass asymp}
\lim_{\varepsilon\to 0}\left[\int_X m\, {\rm d }S+\frac{A_+ A_-}{A}\frac{\log (\varepsilon/4)}{\pi}\right]=\sum_{\pm}\left(\int_{X_\pm}m_\pm \, {\rm d }S+\frac{A_+ A_-}{A}m_\pm (P_\pm)\right).
\end{equation}
In view of \eqref{Zeta Robin}, formula \eqref{int Robin mass asymp} implies \eqref{zeta 1 asymp}.

\subsection*{Justification of asymptotics: Genus 0 example} In the above calculations, we used asymptotic formulas \eqref{abel diff asymp}, \eqref{prime form asymp} proved in \cite{KokTAMS} only for the case in which $x$, $y$ are not close to $Q$, $Q'$. Thus, to justify the above calculations, we need to explain why~\eqref{Phi func asymp} can be integrated over the {\it entire} $X$ including the pinching zone (similarly, why~\eqref{abel diff asymp} can be integrated over the whole contour from $x\in X_\pm$ to $y\in X_\mp$ including the part lying near $Q$). The justification of possibility of these integrations is rather technical and uses some facts from the theory of functions of several complex variables, we preface it with the following simple example illustrating the behaviour of the prime form $E$ for small $\varepsilon$.

Suppose that $X_+=X_-=\overline{\mathbb{C}}$ and $z_\pm$ are standard coordinates on $\mathbb{C}$ (the choice of metrics $\mathfrak{T}_\pm$ is not important but we still assume that $\mathfrak{T}_\pm(z_\pm)=1$ near $z_\pm=0$). Then $X\equiv\overline{\mathbb{C}}$ and $\xi$ is the standard coordinate on $\mathbb{C}$. The prime forms on $X$ and $X_\pm$ are given by
\[\frac{E(\xi,\xi')}{\sqrt{{\rm d}\xi}\sqrt{{\rm d}\xi'}}=\frac{\xi-\xi'}{\sqrt{{\rm d}\xi}\sqrt{{\rm d}\xi'}}, \qquad \frac{E_\pm(z_\pm,z'_\pm)}{\sqrt{{\rm d}z_\pm}\sqrt{{\rm d}z'_\pm}}=\frac{z_\pm-z'_\pm}{\sqrt{{\rm d}z_\pm}\sqrt{{\rm d}z'_\pm}}.\]

If $\xi_1$, $\xi_2$ belong to the same sheet $X'_\pm\subset X$ and $z_1=z_\pm(\xi_1)$, $z_2=z_\pm(\xi_2)$, then formulas \eqref{coordinates near pinch} imply
\begin{align*}
\frac{E(z_1,z_2)^2}{{\rm d}z_1{\rm d}z_2}&=\frac{E(\xi_1,\xi_2)^2}{{\rm d}z_1{\rm d}z_2}\frac{{\rm d}z_1}{{\rm d}\xi_1}\frac{{\rm d}z_2}{{\rm d}\xi_2}=\frac{(\xi_1-\xi_2)^2\xi_1\xi_2}{(1-\xi_1^2)^2(1-\xi_2^2)^2}\frac{4\varepsilon^2}{{\rm d}z_1{\rm d}z_2}.
\end{align*}
In view of the Taylor expansion $\pm\xi=1-\varepsilon/2z_{\pm}+O\bigl(\varepsilon^2/|z_{\pm}|^2\bigr)$, the latter equality can be rewritten as
\begin{align*}
\frac{E(z_1,z_2)^2}{{\rm d}z_1{\rm d}z_2}=&\frac{(z_1-z_2)^2+\cdots}{{\rm d}z{\rm d}z'}=\frac{E_\pm(z_1,z_2)^2+\cdots}{{\rm d}z_1{\rm d}z_2},
\end{align*}
where dots denote the terms growing at most as $O\bigl(\varepsilon\bigl(|z_{1}|^{-1}+|z_{2}|^{-1}\bigr)\bigr)$. Similarly, if $\xi_\pm\in X'_\pm$ and $z_\pm=z_\pm(\xi_\pm)$, then
\begin{align*}
\frac{E(z_+,z_-)^2}{{\rm d}z_+{\rm d}z_-}=-\frac{16z_+^2z_-^2(1+\cdots)}{\varepsilon^2{\rm d}z_+{\rm d}z_-}=-\frac{16E(z_+,0)^2E(z_-,0)^2(1+\cdots)}{\varepsilon^2{\rm d}z_+{\rm d}z_-}.
\end{align*}
Thus, in the genus zero case, asymptotics \eqref{prime form asymp} is valid on $X\times X$ except the zones $|z_\pm(x)|\le\varepsilon^{1-\alpha}$ and $|z_\pm(y)|\le\varepsilon^{1-\alpha}$ where $\alpha>0$ may be arbitrarily small. In this zone, we have
\begin{align*}
\frac{E(z_\pm,z'_{\pm'})^2}{{\rm d}z_\pm {\rm d}z'_{\pm'}}&=\frac{4z_\pm z'_{\pm'}}{\varepsilon^2{\rm d}z_\pm {\rm d}z'_{\pm'}}\Bigl(\sqrt{z'_{\pm'}(z_\pm -\varepsilon)}\mp(\pm')\sqrt{z_\pm (z'_{\pm'}-\varepsilon)}\Bigr)^2,
\end{align*}
whence
\[E(z,z')^2=O\bigl(\varepsilon^{-2\alpha}\bigr), \qquad E(z,z')^{-2}=O\bigl(\varepsilon^2 z^{-2}z'^{-2}(z-z')^{-2}\bigr),\]
and $\log |E(z,z')|^2=-4\pi\Phi(z,z')$ has only logarithmic growth with respect to $z$, $z'$ and $\varepsilon$ while the volumes of the zones $|z_\pm(x)|\le\varepsilon^{1-\alpha}$ and $|z_\pm(y)|\le\varepsilon^{1-\alpha}$ are $O\bigl(\varepsilon^{2(1-\alpha)}\bigr)$. Therefore, one can neglect the contribution of these zones in the integration in \eqref{illegal int}. We show below that $\Phi(z,z')$ behaves similarly near $Q$ in the general case.

\subsection*{Justification of asymptotics}

{\bf 1.} First, we make use of Yamada's theorem (\cite[Lemma~1 and Theorem~1 on p.~116]{Yamada}, see also~\cite[Lemma~1]{KokTAMS}) in the following form.

\begin{Lemma}\label{Yamada lemma}\qquad
\begin{enumerate}\itemsep=0pt
\item[$(i)$] Let $v_\pm$ be holomorphic differentials on $X_\pm$, respectively. Then there is a family $v=v(\varepsilon)$ of holomorphic differentials on $X=X(\varepsilon)$ obeying the inequality
\begin{equation}
\label{Yamada est}
\sum_{\pm}\|v-v_\pm\|_{L_2(X_\pm\backslash U)}+\|v\|_{L_2(U)}\le c\varepsilon,
\end{equation}
where $U=U(\varepsilon)$ is given by
$U:=\{P\in X | |\xi(P)\pm 1|\ge c\}$
and the $L_2$-norm of a one-form~$v$ on subdomain $U\subset X$ is given by
\[\|v\|_{L_2(U)}:=\left(\int_{U}v\wedge\overline{\star v}\right)^{1/2}.\]
\item[$(ii)$] Let $W$ and $W_\pm$ be the canonical bimeromorphic differentials on $X$ and $X_\pm$, respectively. Then the inequality
\begin{equation}
\label{Yamada est merom}
\|W(\cdot,y)-W_\pm(\cdot,y)\|_{L_2(X_\pm\backslash U)}+\|W(\cdot,y)\|_{L_2(X_\mp\cup U)}\le \frac{c\varepsilon}{ d(y,P_\pm)^2}
\end{equation}
holds for any $y\in X_\pm\backslash U$.
\item[$(iii)$] Denote $W_0=(\xi-\xi')^{-2}{\rm d}\xi {\rm d}\xi'$. Then there is a meromorphic differential $Y=Y(\cdot,\xi')$ on~$X$ obeying
\begin{equation}
\label{Yamada est on pinch}
\|Y(\cdot,\xi')\|_{L_2(X\backslash V)}+\|W_0(\cdot,\xi')-Y(\cdot,\xi')\|_{L_2(V)}\le \frac{c}{\big|1-\xi'^2\big|},
\end{equation}
where $\xi'\in V$ and $V$ us given by
$V:=\{P\in X | |z_+(P)|\le c \text{ or } |z_-(P)|\le c\}$.
\end{enumerate}
\end{Lemma}
\begin{proof}
$(i)$ We have
\begin{equation}
\label{Yamada y}
y_{\pm}(P):=\int_{P_\pm}^P v_\pm=\sum_{k>0}c_{k}z_\pm^{k}=\sum_{k>0}\frac{c_{k}\varepsilon^{k}}{\bigl(1-\xi^2\bigr)^k},
\end{equation}
where $z_\pm=z_\pm(P)$, $\xi=\xi(P)$, and $P\in \partial U$.
Let $h=h(\varepsilon)$ be the harmonic extension of $y_\pm|_{\partial U}$ into $U$; then the above expansion implies $\|{\rm d}h\|_{L_{2}(U)}\le c\varepsilon$. Introduce the form $\Phi=\Phi(\varepsilon)$ given by $\Phi=v_\pm$ on $X_\pm\backslash U$ and by $\Phi={\rm d}h$ on $U$. By construction, the jump of the form $\Phi$ on $\partial U$ is zero on any vector field $k$ tangent to $\partial U$. Due to this, ${\rm d}\Phi$ has no singularities on $\partial U$; moreover, we have ${\rm d}\Phi=0$ on the whole $X$. Note that the form $\tilde{\Phi}:=\Phi-{\rm i}\star\Phi$ vanishes outside $U$ due to~${{\rm i}\star {\rm d}z={\rm d}z}$; as a corollary, we have
\begin{equation}
\label{Yamada alpha est h}
\|\tilde{\Phi}\|_{L_2(X)}\le 2\|{\rm d}h\|_{L_2(X)}\le c\varepsilon.
\end{equation}
Consider the Hodge orthogonal decomposition $\tilde{\Phi}={\rm d}\alpha+\delta\beta+\gamma$, where ${\rm d}\gamma=0$, $\delta\gamma=0$ and
\begin{equation}
\label{Yamada alpha est}
\|{\rm d}\alpha\|_{L_2(X)}\le\|\tilde{\Phi}\|_{L_2(X)}\le c\varepsilon.
\end{equation}
Then we have
$\Phi-{\rm d}\alpha=\delta\beta+\gamma+{\rm i}\star\Phi$,
where the left-hand side is closed and the right-hand side is coclosed. So, the form $\Phi-{\rm d}\alpha$ is harmonic and, thus, the form
$v=\frac{1}{2}\bigl(\Phi-{\rm d}\alpha+{\rm i}\star(\Phi-{\rm d}\alpha)\bigr)$
is holomorphic. Note that
\[v-v_\pm=\frac{1}{2}({\rm d}\alpha+{\rm i}\star {\rm d}\alpha) \qquad\text{on}\ X_\pm\backslash U, \qquad v=\frac{1}{2}( d(h-\alpha)+{\rm i}\star d(h-\alpha)) \qquad\text{on}\ U.\]
Then \eqref{Yamada alpha est h}, \eqref{Yamada alpha est} imply \eqref{Yamada est}.

$(ii)$ Repeating the above construction for $v_\pm=W_\pm(\cdot,y)$ and $v_\mp=0$, one gets a differential $v$ obeying \eqref{Yamada est}, where the additional multiplier $ d(y,P_\pm)^{-2}$ should appear in the right-hand side due to the standard estimate $|c_k|\le c|z_\pm(y)|^{-(k+2)}$ for the coefficients in \eqref{Yamada y}. By construction, $v$ and $W(\cdot,y)$ share the same singularity on $X$, i.e., $v-W(\cdot,y)$ is holomorphic on $X$. In addition, the periods of $v$ coincide with those of $(1/2){\rm i}\star {\rm d}\alpha$. Since $\tilde{\Phi}={\rm d}\alpha+\delta\beta+\gamma$ vanishes outside $U$, ${\rm d}\alpha$ is harmonic. Then the increasing smoothness theorems for solutions to elliptic equations imply
\[\Big|\int_l\star {\rm d}\alpha\Big|\le |l|\|{\rm d}\alpha\|_{L_2(X)}\le c|l|\varepsilon d(y,P_\pm)^{-2},\]
where $l$ is any path which does not intersect $V$ and $|l|$ is its length. Due to the last estimate (with $l=a_{k,\pm},b_{k,\pm}$) and the Riemann bilinear relations, we have $v-W(\cdot,y)=\tilde{v}$, where $\|\tilde{v}\|_{L_2(X)}\le c|l|\varepsilon d(y,P_\pm)^{-2}$. As a corollary, we obtain~\eqref{Yamada est merom}.

$(iii)$ To prove $(iii)$, it is sufficient to repeat the above construction with the following replacements: (a)~instead of \eqref{Yamada y}, we consider the anti-derivative $y(P)=-(\xi-\xi')^{-1}$ of $W_0(\cdot,\xi')$ in the domain $V$, (b)~$h$~is defined as the harmonic extension of $y|_{\partial V}$ into domains $X\backslash V=(X_+\backslash V)\cup (X_-\backslash V)$, (c)~$\Phi$~is defined by $\Phi=W_0$ on $V$ and by $\Phi={\rm d}h$ on $X\backslash V$.
\end{proof}

Since the differentials considered in Lemma \ref{Yamada lemma} are holomorphic, formulas \eqref{Yamada est}, \eqref{Yamada est merom}, \eqref{Yamada est on pinch} give rise to {\it point-wise} estimates of $v$, $W$, $Y$. One way to derive them is to combine \eqref{Yamada est}, \eqref{Yamada est merom}, \eqref{Yamada est on pinch} with the mean value property for holomorphic functions,
\begin{equation}
\label{mean value property}
v(z){\rm d}z=\frac{{\rm d}z}{2{\rm i}\pi R^2}\int_{\mathbb{D}_R}{\rm d}\overline{z'}\wedge v(z'){\rm d}z' \Rightarrow |v(z)|\le \frac{\|v\|_{\mathbb{D}_R}}{R\sqrt{2\pi}},
\end{equation}
where $z$, $z'$ are values of the same (arbitrary) holomorphic coordinate and $\mathbb{D}_R$ is the disk ${|z'-z|\le R}$ in $X$.

Put $v_\pm=v_{k,\pm}$, $v_\mp=0$. Due to Lemma \ref{Yamada lemma}\,(i), there is a differential $v$ on $X$ obeying \eqref{Yamada est}. In view of \eqref{mean value property}, the corresponding (i.e., on \smash{$X'_\pm$}) periods of $v$ and $v_{\pm}$ coincide up to the term of order $O(\varepsilon)$ while the remaining periods of $v$ are $O(\varepsilon)$. Therefore, the Riemann bilinear relations imply $\|v-w^\pm_k\|_{L_2(X)}=O(\varepsilon)$.

Now estimate \eqref{Yamada est} yields
\begin{equation}
\label{Yamada est basis}
\|w^\pm_k-v_{k,\pm}\|_{L_2(X_\pm\backslash U)}+\|w^\pm_k\|_{L_2(X_\mp\cup U)}\le c\varepsilon.
\end{equation}
From \eqref{Yamada est basis} and \eqref{mean value property}, it follows that asymptotics \eqref{abel diff asymp} is valid as $\varepsilon/ d(x,P_\pm)\to 0$.

Due to the above fact, the proof of \eqref{abel map asymp} for all $x,y\in X$ is reduced to deriving the estimate
\begin{gather}
\label{abel est}
\max_{x\in U_\alpha}|\mathscr{A}(x-Q)|=o(1),
\end{gather}
where
\[U_\alpha=U_\alpha(\varepsilon):=\{P\in X \mid |z_\pm(P)|\le\varepsilon^{\alpha}\}, \qquad \alpha\in (2/3,1).\]
Since the functions $z_\pm\mapsto v_{k,\pm}(z_\pm)$ are smooth, we have $\|v_{k,\pm}\|_{L_2(U_\alpha\cap X_\pm)}=O\bigl(\varepsilon^{\alpha/2}\bigr)$. Then equation \eqref{Yamada est basis} implies
\begin{equation}
\label{refined yamada est}
\|\vec{\omega}\|_{L_2(U)}\le c\varepsilon, \qquad \|\vec{\omega}\|_{L_2(U_\alpha)}\le c\varepsilon^{\alpha/2}.
\end{equation}
Let us represent $\mathscr{A}(x-Q)$ as the sum of the integrals \smash{$J_1:=\int_{Q}^{x'}\vec{\omega}$} and \smash{$J_2:=\int_{x'}^x\vec{\omega}$}, where the integration path joining $Q$ and $x'$ (resp., $x'$ and $x$) belongs to $U$ (resp., $U_\alpha\backslash U$). Due to \eqref{refined yamada est} and \eqref{mean value property}, we have $|\vec{\omega}(z_\pm)|\le c\varepsilon^{\frac{\alpha}{2}-1}$ on $U_\alpha\backslash U$ and, thus, $|J_2|\le \varepsilon^{\alpha-1}=o(1)$.

 Similarly, from \eqref{refined yamada est} and \eqref{mean value property}, it follows that $|\vec{\omega}(\xi)|\le c\varepsilon$ for $|\xi\pm 1|>c$.

Since $\vec{\omega}$ is holomorphic at $Q$ and $Q'$ (where $\xi(Q)=\infty$ and $\xi(Q')=0$), it obeys $\vec{\omega}(\xi)=O(1)$ near $\xi=0$ and $\vec{\omega}(\xi)=O\bigl(|\xi|^{-2}\bigr)$ near infinity. Thus, the standard estimates of the coefficients in Laurent series for $\vec{\omega}(\xi)$ yield
\[|\vec{\omega}(\xi)|\le\frac{c\varepsilon}{\big|1-\xi^2\big|}, \qquad |\xi\pm 1|>c.\]
Then $|J_1|\le c\varepsilon\int_0^{\infty}\bigl(1+r^2\bigr)^{-1}{\rm d}r=O(\varepsilon)$. Thereby, we have proved \eqref{abel est} and \eqref{abel map asymp}.

Denote
\[\mathfrak{d}(x):=\sqrt{ d(x,Q) d(x,Q')}.\]
By repeating of the above reasoning, we derive the asymptotics
\begin{align*}
W(x,y)=\begin{cases}
W_\pm(x,y)+O(\epsilon), & x,y\in X_\pm, \\
O(\epsilon), & x\in X_\pm, \ y\in X_\mp, \ \epsilon:=\frac{\varepsilon}{\mathfrak{d}(x)\mathfrak{d}(y)}\to 0.
\end{cases}
\end{align*}
from Lemma \ref{Yamada lemma}\,(ii).

In the case of arbitrary $x,y\in X$ (including those close to $Q$ or $Q'$), the above facts, Lemma~\ref{Yamada lemma}\,(iii), and formula \eqref{mean value property} lead to the estimates
\begin{equation}
\label{rought pointwise estimates}
|\vec{\omega}(x)|\le c\varepsilon^{-K}\mathfrak{d}(x)^{-L}, \qquad
|W(x,y)|\le c\varepsilon^{-M}|x-y|^{-2}(\mathfrak{d}(x)\mathfrak{d}(y))^{-L},
\end{equation}
with some {\it finite} $K$. Here $x$, $y$ are values of arbitrary admissible coordinates while the constant~$c$ and the exponents $K$, $M$ in \eqref{rought pointwise estimates} depend on the choice of the coordinate patches of $x$ and $y$; the exact values of $K$, $M$ are of no importance for us.

We make use of the following explicit formula (see \cite[Corollary 2.12]{Fay1})
\begin{equation}
\label{prime form via bidiff}
E^{-2}(x,y)=\frac{\theta^2(e)\bigl(W(x,y)+\sum_{i,j=1}^{g}\partial_i\partial_j\log \theta(e)v_i(x)v_j(y)\bigr)}{\theta\bigl(\int_x^y\vec{v}-e\bigr)\theta\bigl(\int_x^y\vec{v}+e\bigr)}
\end{equation}
for the prime form; here $e\in\mathbb{C}^{g}$ is arbitrary and $\theta=\theta(\cdot|\mathbb{B})$ is the theta function of $X$. Recall that the theta-function $(\vec{z},\mathbb{B}')\mapsto\theta(\vec{z}|\mathbb{B}')$ is smooth in $\vec{z}\in\mathbb{C}^g$ and $\mathbb{B}'$ from the genus $g$ Siegel upper half space. Due to this fact, the equality
\[\theta\bigl((\vec{z}_1,\vec{z}_2)^{\mathsf T}|{\rm diag}(\mathbb{B}_1,\mathbb{B}_2)\bigr)=\theta(\vec{z}_1|\mathbb{B}_1)\theta(\vec{z}_2|\mathbb{B}_2),\]
and the asymptotics \eqref{priod matr asymp}, for each value $x_0$, $y_0$ of admissible coordinates $x$, $y$, there are a~neighbourhood $U\ni(x_0,y_0)$ and a vector $e$ such that the denominator in \eqref{prime form via bidiff} is separated from zero for all sufficiently small $\varepsilon$. Majorizing of the right-hand side of \eqref{prime form via bidiff} with the help of \eqref{rought pointwise estimates} yields the lower bound
\begin{equation}
\label{prime form estimate 1}
c(\log \varepsilon+\log |x-y|+\log \mathfrak{d}(x)+\log \mathfrak{d}(y))\le \log |E(x,y)|,
\end{equation}
 where $c$ is independent of the values $x,y$ but depends on the admissible coordinates they represent.

{\bf 2.} Now we obtain the upper bound for $\log |E(x,y)|$.
First, let us note that the family $\varepsilon\mapsto X_\varepsilon$ of Riemann surfaces is {\it complex analytic} (in the sense of Kodaira \cite[Definition 2.8]{Kod}). Here we consider $\varepsilon$ as a small non-zero {\it complex} parameter, i.e., we allow the rotation of the cuts $\mathscr{C}_\pm(\varepsilon)$ around $P_\pm$ when constructing $X(\varepsilon)$. Due to this fact, the same reasoning as in \cite[Section~2.3]{KokKorMPAG} or \cite[Section~2]{KokKorJDG} (see also \cite[Chapter~3]{Fay}) shows that all the canonical objects on $X=X(\varepsilon)$ (such as basic Abelian differentials $w=w_{k}^\pm$, the canonical bimeromorphic differential $W$, the square of the prime form~$E^2$, etc.) are complex analytic with respect to $\varepsilon\ne 0$ and the admissible coordinates. Thus, we can use the results of the theory of functions of several complex variables (see, e.g.,~\cite{BM}) in deriving the asymptotics of the prime form~$E$ as $\varepsilon\to 0$.

Due to \cite[Theorem 2]{KokTAMS}, the asymptotics
\begin{equation}
\label{merom diff asymp}
W(x,y)=\begin{cases}
W_\pm(x,y)+\dfrac{\varepsilon^2}{16}W_\pm(x,z_\pm)W_\pm(y,z_\pm)|_{z_\pm=0}+o\bigl(\varepsilon^2\bigr), & x\in X_\pm, \ y\in X_\pm,\\
-\dfrac{\varepsilon^2}{16}W_\pm(x,z_\pm)W_\mp(y,z_\mp)|_{z_\pm=z_\mp=0}+o\bigl(\varepsilon^2\bigr), & x\in X_\pm, \ y\in X_\mp,
\end{cases}
\end{equation}
as $\varepsilon\to 0$ is valid for $x$, $y$ lying outside arbitrarily small fixed neighbourhood of $Q$. Formulas \eqref{abel diff asymp}, \eqref{merom diff asymp} and the Hartogs's theorem on separate holomorphicity imply that $w(z_\pm)$ and $W_\pm(z_\pm,z'_\pm)$ are complex analytic in $z_\pm$, $z'_\pm$, $\varepsilon$ in the domains
\begin{align*}
\mathscr{D}_1(z_\pm):&=\{(z_\pm,\varepsilon) \mid c_0|\varepsilon|<|z_\pm|<c_1\},\\
\mathscr{D}_2(z_\pm):&=\{(z_\pm,z'_\pm,\varepsilon) \mid c_0|\varepsilon|<|z_\pm|,\,|z'_\pm|<c_1, \, z_\pm\ne z'_\pm\},
\end{align*}
respectively. Since $\mathscr{D}_{1,2}(z_\pm)$ are Hartogs domains, $w(z_\pm)$ and $W(z_\pm,z'_\pm)$ admit expansions in the convergent Hartogs series
\begin{equation}
\label{Hartogs series z}
w(z_\pm)=\sum_{k\ge 0}\varepsilon^k w^{(k)}(z_\pm), \qquad W(z_\pm,z'_\pm)=\sum_{k\ge 0}\varepsilon^{k} W^{(k)}(z_\pm,z'_\pm),
\end{equation}
where the coefficients $w^{(k)}$ and $W^{(k)}$ are independent of $\varepsilon$ and holomorphic in the base domains
\[\{0<|z_\pm|<c_1\}, \qquad \{0<|z_\pm|,|z'_\pm|<c_1,\, z_\pm\ne z'_\pm\},\]
respectively.

Now we are to study the behaviour of the coefficients in \eqref{Hartogs series z} near singularities $z_\pm=0$ and~${z'_\pm=0}$. To this end, we make use of the Cauchy formula
\begin{equation}
\label{Cauchy formula}
w^{(l)}(z_\pm)=\frac{1}{2\pi {\rm i} }\oint\limits_{|\varepsilon|=\varrho_0}\frac{w(z_\pm)\,{\rm d}\varepsilon}{\varepsilon^{l+1}}, \qquad W^{(l)}(z_\pm,z'_\pm)=\frac{1}{2\pi {\rm i} }\oint\limits_{|\varepsilon|=\varrho_0}\frac{W(z_\pm,z'_\pm)\,{\rm d}\varepsilon}{\varepsilon^{l+1}},
\end{equation}
where $\varrho_0>0$ is sufficiently small. Note that $W(z_\pm,z'_\pm)$ is well defined for $|z_\pm|,|z'_\pm|\ge c_0|\varepsilon|$. We put $|z_\pm|=c|\varepsilon|$ or $|z'_\pm|=c|\varepsilon|$ into \eqref{Cauchy formula} and majorize the right-hand sides by means of \eqref{rought pointwise estimates}. As a result, we obtain
\begin{gather*}
\bigl|w^{(l)}(z_\pm)\bigr|=O\bigl(|z_\pm|^{-(K+L/2+l)}\bigr),\qquad
\bigl|W^{(l)}(z_\pm,z'_\pm)\bigr|=O\bigl(|z_\pm z'_\pm|^{-(M+L/2+l)}|z_\pm-z'_\pm|^{-2}\bigr),
\end{gather*}
i.e., the coefficients in \eqref{Cauchy formula} have no essential singularities at $z_\pm=0$ or $z'_\pm=0$.

Rewrite \eqref{Hartogs series z} as
\begin{gather}
w(z_\pm)=\frac{1}{z_\pm^{M}}\sum_{k\ge 0}\left(\frac{\varepsilon}{z_\pm}\right)^{k}\tilde{w}^{(k)}(z_\pm), \nonumber\\
W(z_\pm,z'_\pm)=\frac{1}{(z_\pm-z'_\pm)^{2}(z_\pm z'_\pm)^{M}}\sum_{k\ge 0}\left(\frac{\varepsilon}{z_\pm z_\pm'}\right)^{k}\tilde{W}^{(k)}(z_\pm,z'_\pm),\label{Hartogs series z new}
\end{gather}
where the coefficients $\tilde{w}^{(k)}(z_\pm)$ and $\tilde{W}^{(k)}(z_\pm,z'_\pm)$ are holomorphic at $z_\pm=0$ or $z'_\pm=0$. Substitution of series \eqref{Hartogs series z new} into \eqref{prime form asymp} yields the expansion
\[E^{2}(z_\pm,z'_\pm)=\varepsilon^n(z_\pm z'_\pm)^N(z_\pm-z'_\pm)^2\left(\mathscr{E}(z_\pm,z'_\pm)+O\left(\frac{\varepsilon}{z_\pm z'_\pm}\right)\right),\]
where $|z_\pm z'_\pm|\ge c|\varepsilon|\to 0$ and $n$, $N$ and $\mathscr{E}$ are independent of $\varepsilon$. As a corollary, we obtain the estimate
\begin{equation}
\label{prime form estimate 2}
\log |E^2(x,y)|\le C\left(|{\log \varepsilon}|+|{\log |z_\pm(x)}|+|{\log |z_\pm(y)}|+\log |z_\pm(x)-z_\pm(y)|\right)
\end{equation}
for the case $|z_\pm z'_\pm|\ge c|\varepsilon|$.

To prove \eqref{prime form estimate 2} for the case $|z_\pm z'_\pm|< c|\varepsilon|$, it is sufficient to repeat the above reasoning for the following Hartogs series and domains
\begin{gather*}
w(\zeta)=\sum_{k\ge 0}\varepsilon^k w^{(k)}(\zeta), \qquad \bigl(c_0|\varepsilon|<|\zeta|<c_1|\varepsilon|^{-1}\bigr),\\
W(\zeta,\zeta')=\sum_{k\ge 0}\varepsilon^{k} W^{(k)}(\zeta,\zeta'),\qquad \bigl(c_0|\varepsilon|<|\zeta|, |\zeta'|<c_1|\varepsilon|^{-1},\, \zeta\ne\zeta'\bigr),\\
W(\zeta,z'_\pm)=\sum_{k\ge 0}\varepsilon^{k} W^{(k)}(\zeta,z'_\pm), \qquad \bigl(c_0|\varepsilon|<|\zeta|<c_1|\varepsilon|^{-1},\, c_0|\varepsilon|<|z_\pm'|<c_1\bigr).
\end{gather*}

Now combining estimates \eqref{prime form estimate 1} and \eqref{prime form estimate 2} and formula \eqref{Phi func} leads to
\[|\Phi(x,y)|\le C\left(|\log \varepsilon|+|\log d(x,y)|+|\log d(x,Q)|+|\log d(y,Q)|\right).\]
Due to this estimate, for $\alpha\in (2/3,1)$ the integral of $\Phi(x,\cdot)$ over a $\varepsilon^{\alpha}$-neighbourhood, $U_\alpha$, of $Q$ (in the metric $\mathfrak{T}$) is $O\bigl(\varepsilon^{3\alpha/2}\bigr)$. Hence, one can neglect the integration over $U_\alpha\times X$ or $X\times U_\alpha$ in~\eqref{illegal int}. By this, we justify the calculations leading to \eqref{zeta 1 asymp}. Thus, we have proved \eqref{zeta 1 asymp}.

\subsection*{Asymptotics of the first non-zero eigenvalue $\boldsymbol{\lambda_1(\varepsilon)}$ of $\boldsymbol{\Delta(\varepsilon)}$} Consider the equation
\begin{equation}
\label{first eigenpair}
(\Delta-\lambda_1)u_1=0
\end{equation}
in $X$. To describe the asymptotics of solution to \eqref{first eigenpair} as $\varepsilon\to 0$, we apply the method of matched expansions (in the form given in \cite[Chapter 6]{MNP}).

Far away from $Q$, we use the solution $U_\pm:=c_\pm R_{\lambda_1,\pm}(z_\pm,P_\pm)$ to \eqref{first eigenpair} in $X_\pm$ as an approximation for $u_1$ (where $R_{\cdot,\pm}$ are the resolvent kernel for $\Delta_\pm$). Note that $U_\pm$ grows logarithmically near $Q$.

Near $Q$ we introduce the new coordinates
\smash{$q_\pm=\frac{2z_\pm}{\varepsilon}-1$}, \smash{$ k=q_\pm\pm\sqrt{q_\pm^2-1}$}
(the latter is the inverse Joukowsky transform of the $q_\pm$). In coordinates $q_\pm$ equation \eqref{first eigenpair} takes the form~${\bigl(\partial_{q_\pm}\partial_{\overline{q_{\pm}}}-\lambda_1\varepsilon^{2}/4\bigr)=0}$.
Thus, the approximation $w(q_\pm)$ of $u_1$ near $Q$ should be harmonic. In order to make $w$ non-trivial, we allow it to grow logarithmically as $q_\pm\to\infty$. Then the required $w$ is of the form
\[w=D\log |k|+B=D\log \big|q_\pm\pm\sqrt{q_\pm^2-1}\big|+B.\]

Now, we should ``glue'' the above approximations together in the intermediate zone
\[|z_\pm|\asymp\varepsilon^{1/2} \ \Leftrightarrow \ |q_\pm|\asymp \varepsilon^{-1/2}.\]
To this end, the asymptotics of approximations as $z_\pm\to 0$ and $q_\pm\mapsto \infty$ should be consistent in the intermediate zone. From the expansion $R_{\lambda,\pm}=-\frac{1}{A_\pm\lambda}+G_\pm+O(\lambda)$, it follows that
\begin{align}
U_\pm&=-\frac{c_\pm}{A_\pm\lambda_1}\bigl(1-A_\pm\lambda_1 G(z_\pm,P_\pm)+O\bigl(\lambda_1^2\bigr)\bigr)\nonumber\\
&=-\frac{c_\pm}{A_\pm\lambda_1}\left(1+\frac{A_\pm\lambda_1}{2\pi} \log |z_\pm|-A_\pm\lambda_1 m_\pm(P_\pm)+O\bigl(\lambda_1^2\bigr)+o(1)\right)\label{asump of ext sol near int}
\end{align}
near $Q$. In view of the Taylor expansion
\[\sqrt{q_\pm^2-1}=\frac{2z_\pm}{\varepsilon}\sqrt{1-\frac{\varepsilon}{z_\pm}}=\frac{2z_\pm}{\varepsilon}-1-\frac{\varepsilon}{4z_\pm}+O\left(\Big|\frac{\varepsilon}{z_\pm}\Big|^2\right),\]
the asymptotics of $w(q_\pm)$ as $\varepsilon/z_\pm\asymp\varepsilon^{1/2}\to 0$ is given by
\begin{equation}
\label{asump of int sol near ext}
w(q_\pm)=\pm D\bigl(\log |z_\pm|-\log (\varepsilon/4)\bigr)+B+O\bigl(\varepsilon^{1/2}\bigr) .
\end{equation}
{\samepage Matching asymptotics \eqref{asump of ext sol near int} and \eqref{asump of int sol near ext} yields
\[\mp 2\pi D=c_\pm=\frac{B\mp D\log (\varepsilon/4)}{-(A_\pm\lambda_1)^{-1}+m_\pm(P_\pm)+O(\lambda_1)},\]
whence
\[-\frac{1}{\lambda_1 A_\pm}+O(\lambda_1)+m_\pm(P_\pm)=\frac{1}{2\pi}\log (\varepsilon/4)\mp \frac{B}{2\pi D}.\]
Now, the summation over $\pm$ yields the asymptotics
\[-\frac{1}{\lambda_1}\left(\frac{A}{A_+A_-}\right)+O(\lambda_1)=\frac{1}{\pi}\log (\varepsilon/4)-\sum_\pm m_\pm(P_\pm)\]
leading to \eqref{first eigenvalue asym}.}

Then the global approximate solution to \eqref{first eigenpair} is constructed by multiplying the (matched) local approximations $U_\pm$, $w$ by the appropriate smooth cut-off functions and adding them together. We do not describe this (well known) construction in detail since it can be found in \cite[Chapter~6]{MNP} (or, e.g., in \cite{KorSIF}). As a result, we obtain the function $\hat{u}_1$ obeying \smash{$\big\|\big(\Delta-\hat{\lambda}_1\big)\hat{u}_1\big\|_{L_2(X)}=o(1)$} and \smash{$\|\hat{u}_1\|_{L_2(X)}=1+o(1)$}, where $\hat{\lambda}_1$ obeys \eqref{first eigenvalue asym} with omitted term $o(1)$. Then the same reasoning as in \cite[Section 4.2]{MNP} (see also \cite{KorMax}) leads to estimate \smash{$\big|\lambda_1-\hat{\lambda}_1\big|+\big\|u_1-\hat{u}_1\big\|_{L_2(X)}=o(1)$}. By this, we complete the proof of \eqref{first eigenvalue asym}. The same arguments (with $\lambda_1$ replaced by $\lambda_k$, $k>1$) lead to \eqref{separated eigenvalues asym}.

\subsection*{Pointwise asymptotics for Green function} For simplicity, suppose that $x\in X_+$ is separated from $P_+$. Integration of \eqref{Green Verlinde} yields
\begin{align}
\label{scalar Robin}
m(x)+\langle m\rangle=-\frac{2}{A}\int_{X}\Phi(x,y)\, {\rm d }S_\rho(y),
\end{align}
where $\langle m\rangle:=\frac{1}{A}\int_X m{\rm d }S$ is the average Robin mass; note that \eqref{scalar Robin} remains valid after the replacement $m,X,A,\Phi\to m_\pm,X_\pm,A_\pm,\Phi_\pm$. In view of \eqref{scalar Robin},\eqref{Phi func asymp}, and \eqref{int Robin mass asymp}, we have
\begin{gather*}
Am(x)+\sum_{\pm}\left(A_\pm\langle m_\pm\rangle+\frac{A_+A_-}{A} m_\pm(P_\pm)\right)-\frac{A_+A_-}{A}\frac{\log (\varepsilon/4)}{\pi}\\
\qquad=o(1)
-2\int_{X_+}\Phi_+(x,\cdot)\,{\rm d }S-2\int_{X_-}\left(\Phi_+(x,P_+)+\Phi_-(y,P_-)+\frac{\log (\varepsilon/4)}{2\pi}\right){\rm d }S(y)\\
\qquad=o(1)
+\sum_{\pm}A_\pm\langle m_\pm\rangle+A_+m_+(x)-2A_-\Phi_+(x,P_+)+A_-m_-(P_-)-\frac{AA_-}{A}\frac{\log (\varepsilon/4)}{\pi},
\end{gather*}
whence
\begin{align}
m(x)={}&\frac{A_+m_+(x)-2A_-\Phi_+(x,P_+)}{A}-\frac{A_+A_-}{A^2}m_+(P_+)\nonumber\\
&+\frac{A_-^2}{A^2}\left(m_-(P_-)-\frac{\log (\varepsilon/4)}{\pi}\right)+o(1).\label{pointwise robin asymp}
\end{align}

Suppose that $x,y\in X_+$ are separated from $P_+$ and put $Q=(m_+(x)+m_+(y))/2$. Then \eqref{Green Verlinde}, \eqref{Phi func asymp} and \eqref{pointwise robin asymp} imply
\begin{align*}
G(x,y)={}&\Phi_+(x,y)+\frac{A_+Q}{A}-\frac{A_-}{A}(\Phi_+(x,O_+)+\Phi_+(y,O_+))-\frac{A_+A_-}{A^2}m_+(P_+)\\
&+\frac{A_-^2}{A^2}\left(m_-(P_-)-\frac{\log (\varepsilon/4)}{\pi}\right)+o(1)\\
={}&\frac{A_-^2}{A^2}\left(\sum_\pm m_\pm(P_pm)-\frac{\log (\varepsilon/4)}{\pi}\right)
+G_+(x,y)\\
&-\frac{A_-}{A}(G_+(x,O_+)+G_+(y,O_+))+o(1).
\end{align*}
Similarly, if $x\in X_\pm$, $y\in X_\mp$ are separated from $P_\pm$, $P_\mp$, respectively, we obtain
\begin{align*}
\lim_{t\to 0}\left[G(x,y)+\frac{A_+A_-}{A^2}\left(\sum_{\pm}m_\pm(P_\pm)-\frac{\log (\varepsilon/4)}{\pi}\right)\right]=\sum_{\pm}\frac{A_\pm G_\pm(x,O_\pm)}{A}+o(1).
\end{align*}

\begin{Remark} It is interesting to compare the results of this Section with those from \cite{Went} and~\cite{JWent}. Although the scheme of degeneration of the Riemann surface in \cite{JWent,Went} differs from ours (we use cross-gluing along the straight cuts, whereas, say in \cite{Went}, the standard pinching construction from \cite{Fay1} is employed) and the metrics (e.g., Bergman's, Arakelov's) considered in these papers differ significantly from singular Troyanov's, a certain similarity can be observed. In particular, formula~\eqref{separated eigenvalues asym} may be compared with ``spectral convergence'' from \cite{JWent}. We also notice that the author of \cite{Went} had to overcome the same difficulty (finding estimates for the integrals over the pinching zone) as we were struggling with in the main part of the present section.
\end{Remark}

\begin{Remark}\label{smoth conical case rem 3}
All of the above constructions and arguments remain valid if the metrics on $X_\pm$ are smooth outside conical singularities and flat near them and the points $P_\pm$, respectively. In~view of Remarks~\ref{smoth conical case rem 1} and \ref{smoth conical case rem 2}, this means that asymptotics \eqref{zeta 1 asymp}--\eqref{separated eigenvalues asym} remain valid in this case.
\end{Remark}

\begin{Remark}
There are the following probabilistic interpretations of the pointwise Robin mass~$m(y)$ and $\operatorname{reg}\zeta_\Delta(1)$ (see \cite[Theorem 4.2]{Steiners-Game}). Suppose that a Brownian particle starts from some given point $x\in X$, then the expectation value $t(x\to {\rm rand})$ of the time it takes to approach the (metric) $\epsilon$-neighbourhood of a randomly chosen point $y\in X$ is independent of $x$ and equals~to
\[t(x\to {\rm rand})=\frac{\log (2/\epsilon)-\gamma}{2\pi}A+\operatorname{reg}\zeta_\Delta(1)+O(\epsilon).\]
At the same time, if a Brownian particle starts from a random point of $X$, then the expectation value $t({\rm rand}\to y)$ of the time it takes to approach the (metric) $\epsilon$-neighbourhood of a given point~${y\in X}$ is equal~to
\[t({\rm rand}\to y)=-\frac{A\log \epsilon}{2\pi}+m(y)+O(\epsilon).\]
Therefore, the logarithmic growth of $m(y)$ and $\operatorname{reg}\zeta_\Delta(1)$ as $\varepsilon\to 0$ provided by \eqref{zeta 1 asymp} and \eqref{pointwise robin asymp} can be interpreted as the {\rm(}logarithmic{\rm)} growth of the average time it takes for a Brownian particle to move across the pinching zone (from $X'_\pm$ to $X'_\mp$).
\end{Remark}

\section[Computing of reg zeta(1) for self-adjoint extensions of symmetric Laplacian on a tetrahedron with conical angles pi]{Computing of $\boldsymbol{\operatorname{reg}\zeta(1)}$ for self-adjoint extensions of\\ symmetric Laplacian on a tetrahedron with conical angles $\boldsymbol{\pi}$}
\label{Examples sec}
It is instructive to consider a simple example of a genus zero polyhedral surface with four conical points of conical angles $\pi$ (a tetrahedron). Due to availability of explicit information about the spectrum, this example can be completely analysed by means of a direct computation independent of the previous considerations. Moreover, we also introduce the exotic self-adjoint extensions of the symmetric Laplacian on the tetrahedron that are non-Friedrichs and compute the $\operatorname{reg}\zeta(1)$ for them.

 Let $\mathbb{T}=\mathbb{C}/\Lambda$ ($\Lambda=\mathbb{Z}+\tau\mathbb{Z}$) be a torus with periods $1,\tau$, endowed with the standard flat metrics (then its area is equal to $\operatorname{Im}\tau$). Let $AA'A'''A''$ ($A=0$, $A'=1$, $A''=\tau$, $A'''=(1+\tau)$) be the fundamental parallelogram of $\mathbb{T}$.

Consider the holomorphic involution, $\dag$, of $\mathbb{T}$
$\dag\colon z \,{\rm mod}\, \Lambda\mapsto -z \,{\rm mod}\, \Lambda$.
Then $\pi\colon \zeta\mapsto\zeta/\dag$ is the isometric double cover from $\mathbb{T}$ onto $\mathbb{T}/\dag$ ramified at $A$ and $B=1/2$, $C=\tau/2$ and $D=(1+\tau)/2$ (mod $\Lambda$).

The base $\mathbb{T}/\dag$ coincides with the tetrahedron $\mathscr{T}=ABCD$ obtained from its net, the triangle $AA'A''$, by gluing of triangles $ABC$, $BCD$, $BA'D$ and $A''CD$ along equal sides. Each tetrahedron (after suitable homothety) with all conical angles $\pi$ can be obtained in such a way.

\begin{figure}[!ht]
\center{\includegraphics[width=0.5\linewidth]{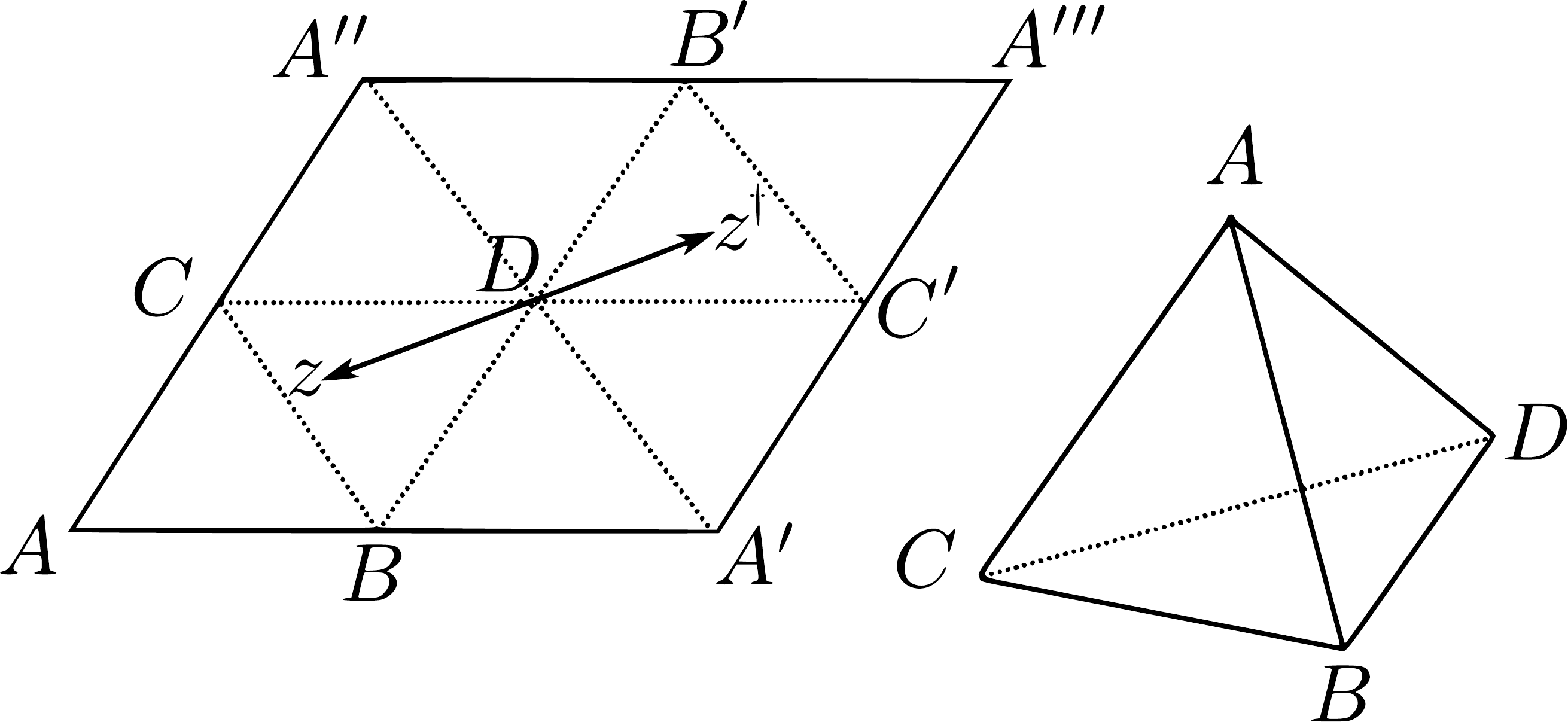}}
\caption{The tetrahedron $\mathscr{T}$ and its double cover $\mathbb{T}$.}\label{fignya}
\end{figure}

Let $\Delta_\mathbb{T}$ and $\Delta_4$ be the Laplacians on $\mathbb{T}$ and $\mathscr{T}$, respectively. Obviously, $(\lambda,u)$ is an eigenpair of $\Delta_4$ if and only if $(\lambda,u\circ\pi)$ is an eigenpair of $\Delta_\mathbb{T}$.

 Conversely, let us search for eigenfunctions of $\Delta_\mathbb{T}$ in the form $v={\rm e}^{{\rm i}(az+b\overline{z})}$; then $\lambda=4ab$ and the periodicity of $v$ with respect to the lattice $\Lambda$ implies $a+b=2\pi n$, $a\tau+b\overline{\tau}=2\pi m$. Since the linear space spanned by such eigenfunctions is an algebra closed under complex conjugation, it is dense in $C(T)$ and $L_2(T)$ (this means that all eigenfunctions of $\Delta_\mathbb{T}$ are of the above form). As a corollary, each non-zero eigenvalue $\lambda=4ab$ of $\Delta_\mathbb{T}$ is double and its eigenspace is spanned by the eigenfunctions $\cos (az+b\overline{z})$, $\sin(az+b\overline{z})$ which are (respectively) even and odd under the involution $v\mapsto v^*=v\circ\dag$. In particular, the non-zero eigenvalues of $\Delta_\mathbb{T}$ are eigenvalues of $\Delta_4$ with double multiplicity.

 Thus, we have
$%\label{tor-sphere-main}
2\zeta_{\Delta_4}\equiv\zeta_{\Delta_\mathbb{T}}$ and
\begin{equation}
\label{doubling zetas}
\operatorname{reg}\zeta_{\Delta_4}(1)=\frac{1}{2}\operatorname{reg}\zeta_{\Delta_\mathbb{T}}(1).
\end{equation}

To calculate $\operatorname{reg}\zeta_{\Delta_\mathbb{T}}(1)$, one can use the explicit formula for the Green function of $\Delta_\mathbb{T}$ on \cite[p.~21]{Fay} which leads to the expression
\begin{equation}
\label{Robin mass torus}
m(\mathbb{T})=-\frac{\log |2\pi\eta(\tau)|}{\pi}
\end{equation}
for the Robin mass on $\mathbb{T}$ (here $\eta$ is the Dedekind eta function). Combining this formula with~\eqref{Zeta Robin} and \eqref{doubling zetas}, one arrives at
\[%\label{Okio}
\operatorname{reg}\zeta_{\Delta_4}(1)=\frac{A_4}{\pi}\left[\frac{\gamma-\log 2}{2}-\log |2\pi\eta(\tau)|\right],
\]
where $A_4=\operatorname{Im}\tau/2$ is the area of the tetrahedron. Thus, considering tetrahedron of unit area with aspect ratio \eqref{aspect ratio}, we arrive at~\eqref{zeta tetra}.

\begin{Remark}
The explicit expression for $\operatorname{reg}\zeta_{\Delta_\mathbb{T}}(1)$ was first obtained in \cite[equations (A.1)--(A.4), p.~800]{Okikitorus}.
\end{Remark}

\begin{Remark} In \cite{Okikitorus}, it was stated that, for a given area, $\operatorname{reg}\zeta_{\Delta_\mathbb{T}}(1)$ attains (global) minimum at the hexagonal torus \smash{$\tau=\frac{1+{\rm i}\sqrt{3}}{2}$}. The latter statement essentially coincides with claim of Osgood, Phillips and Sarnak (see \cite[Lemma~4.1 and a numerical evidence given on p.~206]{OPS1})
concerning the global maximum \big(at \smash{$\tau=\frac{1+{\rm i}\sqrt{3}}{2}$}\big) of the function $f(\tau)= \operatorname{Im} \tau |\eta(\tau)|^4$ on the fundamental domain of the modular group. It is interesting to notice that, although the graph of the function~$f$ (which can be obtained, using a modern laptop, in a minute, see, e.g., \cite[Figure~2, p.~78]{KleinKokKor}) shows the spectacular global maximum at \smash{$\frac{1+{\rm i}\sqrt{3}}{2}$}, the formal mathematical proof of the latter {\it global} maximality, to the best of our knowledge, is still lacking.
\end{Remark}

\begin{Remark}
The global maximality at \smash{$\frac{1+{\rm i}\sqrt{3}}{2}$} of the function $f$ from the previous remark would imply that among tetrahedrons with conical angles $\pi$ and with given area, the (global) minimum of $\operatorname{reg}\zeta_{\Delta_4}(1)$ is attained at the regular tetrahedron.
\end{Remark}

\subsection*{Calculation of $\boldsymbol{\operatorname{reg}\zeta(1)}$ for self-adjoint extensions of Laplacian on $\boldsymbol{\mathscr{T}}$}
Let $\dot{\Delta}_4$ be the $L_2(\mathscr{T})$-closure of the Laplacian defined on the set $\dot{C}^{\infty}(\mathscr{T})$ of smooth functions on the tetrahedron $\mathscr{T}\simeq\overline{\mathbb{C}}$ vanishing near one of the vertices (say, $A$). Since all conical angles of $\mathscr{T}$ are less than $2\pi$, all the self-adjoint extensions of $\dot{\Delta}_4$ are the Laplacians $\Delta_{4,\alpha}$ ($|\alpha|\le\pi/2$) with the domains spanned by elements of $\operatorname{Dom}\dot{\Delta}_4$ and the function
\begin{equation}
\label{domain pseudo}
\chi\left(\frac{\log r}{2\pi}\sin(\alpha)+\cos (\alpha)\right),
\end{equation}
on $\mathscr{T}$, where $r$ is the geodesic distance to the vertex and $\chi$ is a smooth cut-off function with sufficiently small support equal to one near $A$ \cite{HK,KL}.

 Due to the divergence theorem, ${\rm Ker}\Delta_\alpha=\{0\}$ except the case $\alpha=0$ (corresponding to the Friedrichs extension $\Delta_{4,0}:=\Delta_4$).

Now, let $\dot{\Delta}_{\mathbb{T}}$ be the $L_2(\mathscr{T})$-closure of the Laplacian defined on the set $\dot{C}^{\infty}(\mathbb{T})$ of smooth functions vanishing near $A=0 ({\rm mod}\Lambda)$. Then all its self-adjoint extensions are the pseudo-Laplacians (see \cite{CdV}) $\Delta_{\mathbb{T},\alpha}$ ($|\alpha|\le\pi/2$) with domains spanned by elements of $\operatorname{Dom}\dot{\Delta}_{\mathbb{T}}$ and function~\eqref{domain pseudo} on~$\mathbb{T}$.

Clearly, $u\in\operatorname{Dom}\Delta_{4,\alpha}$ if and only if $u\circ\pi\in\operatorname{Dom}\Delta_{\mathbb{T},\alpha}$.

 At the same time, any odd function $v\in\operatorname{Dom}\Delta_{\mathbb{T},\alpha}$ has zeroes at all the ramification points~$A$, $B$, $C$, $D$ (modulo $\Lambda$) and, therefore, belongs to the domain of any pseudo-Laplacian $\Delta_{\mathbb{T},\alpha'}$.

Therefore, the eigenvalues of $\Delta_{\mathbb{T},\alpha}$ can be divided into two series
\begin{enumerate}\itemsep=0pt
\item[(1)] the eigenvalues with even eigenfunctions which coincide with the eigenvalues of $\Delta_{4,\alpha}$,
\item[(2)] the eigenvalues with odd eigenfunctions which coincide with the eigenvalues of $\Delta_{4}$.
\end{enumerate}

Notice that the above division of the spectrum of $\Delta_{\mathbb{T},\alpha}$ agrees with \cite[Theorem 2]{CdV}, presenting an illustration to it.

As a corollary, we have $\zeta_{\Delta_{\mathbb{T},\alpha}}\equiv\zeta_{\Delta_{4,\alpha}}+\zeta_{\Delta_4}$ and
\begin{equation}
\label{zetas pseudo}
\operatorname{reg}\zeta_{\Delta_{\mathbb{T},\alpha}}(1)=\operatorname{reg}\zeta_{\Delta_{4,\alpha}}(1)+\operatorname{reg}\zeta_{\Delta_{4}}(1).
\end{equation}
Here the (regularized) $\zeta$-functions are defined by formulas \eqref{zeta}, \eqref{reg zeta of 1}, where $\Delta$ and its nonzero eigenvalues is replaced by $\Delta_{\mathbb{T},\alpha}$, $\Delta_{4,\alpha}$, or $\Delta_4$, and their nonzero eigenvalues, respectively. Note that the values $\operatorname{reg}\zeta_{\Delta_{\mathbb{T},\alpha}}(1)$, $\operatorname{reg}\zeta_{\Delta_{4,\alpha}}(1)$ ($\alpha\ne 0$) are independent of the choice of the branch $\lambda_1^{s}=|\lambda|^s {\rm e}^{\pm\pi {\rm i} s}$ of $\lambda_1^{s}$, where $\lambda_1<0$ is the first eigenvalue of both $\Delta_{\mathbb{T},\alpha}$ and $\Delta_{4,\alpha}$.

The left-hand side of \eqref{zetas pseudo} can be calculated by the use of the following formula of \cite{AisHK}:
\begin{equation}
\label{comparison aissi}
{\rm Tr}\bigl((\Delta_{\mathbb{T},\alpha}-\lambda)^{-1}-(\Delta_{\mathbb{T}}-\lambda)^{-1}\bigr)=\partial_\lambda\log \bigl(m_\lambda+\operatorname{cot}\alpha\bigr),
\end{equation}
where $\lambda$ does not belong to the spectra of $\Delta_{\mathbb{T},\alpha}$ and $\Delta_{\mathbb{T}}$, and $m_\lambda(\mathbb{T})$ is the coefficient in the asymptotics of the resolvent kernel $R_\lambda$ of $\Delta_{\mathbb{T}}$ near $A$
\[R_\lambda(x,A)=-\frac{\log r}{2\pi}+m_\lambda(\mathbb{T})+o(1).\]
Since zero is an eigenvalue of $\Delta_T$ but not of $\Delta_{\mathbb{T},\alpha}$, we have
\[
%\label{formal reg of zeta in lambda}
\operatorname{reg}\zeta_{\Delta_{\mathbb{T},\alpha}}(1)-\operatorname{reg}\zeta_{\Delta_{\mathbb{T}}}(1)=\lim_{\lambda\to 0}\left[{\rm Tr}\bigl((\Delta_{\mathbb{T},\alpha}-\lambda)^{-1}-(\Delta_{\mathbb{T}}-\lambda)^{-1}\bigr)-\frac{1}{\lambda}\right].
\]
At the same time, from the expansion $R_\lambda=-\frac{1}{A_{\mathbb{T}}\lambda}+G+O(\lambda)$ for the resolvent kernel as $\lambda\to 0$, we obtain
\[m_\lambda(\mathbb{T})=-\frac{1}{A_{\mathbb{T}}\lambda}+m(\mathbb{T})+O(\lambda),\]
where $A_{\mathbb{T}}$ is the area of $\mathbb{T}$ and the Robin mass $m(\mathbb{T})$ is given by \eqref{Robin mass torus}. Then the right-hand side of \eqref{comparison aissi} is
\[\frac{1}{\lambda}-A_{\mathbb{T}}(\operatorname{cot}\alpha+m(\mathbb{T}))\lambda+O\bigl(\lambda^2\bigr).\]
Then we obtain
\begin{equation}
\label{comp form reg zetas torus}
\operatorname{reg}\zeta_{\Delta_{\mathbb{T},\alpha}}(1)-\operatorname{reg}\zeta_{\Delta_{\mathbb{T}}}(1)=-A_{\mathbb{T}}(\operatorname{cot}\alpha+m(\mathbb{T})).
\end{equation}
Combining \eqref{zetas pseudo} with \eqref{comp form reg zetas torus} and \eqref{doubling zetas}, one arrives at
\begin{align*}
\operatorname{reg}\zeta_{\Delta_{4,\alpha}}(1)=\operatorname{reg}\zeta_{\Delta_{4}}(1)-2A_{4}(\operatorname{cot}\alpha+m(\mathbb{T})).
\end{align*}
Now, considering the tetrahedron of unit area with aspect ratio \eqref{aspect ratio}, we obtain \eqref{reg zeta for pseudoL on tetra} from~\eqref{comp form reg zetas torus} and \eqref{zeta tetra}, \eqref{Robin mass torus}.

\subsection*{Acknowledgements}
The authors thank the anonymous referees for careful reading of the manuscript, valuable comments and significant improvements.
The first author thanks Max Planck Institute for Mathematics in Bonn for warm hospitality and excellent working conditions. The research of the second author was supported by Fonds de recherche du Qu\'ebec.

\pdfbookmark[1]{References}{ref}
\LastPageEnding

\end{document}